\documentstyle[12pt]{amsart}
\newsymbol\subsetneq 2328
\newsymbol\subsetneqq 2324
\newsymbol\supsetneq 2329
\newsymbol\nsubseteq 232A

\newtheorem{prop}{Proposition}[section]
\newtheorem{th}[prop]{Theorem}
\newtheorem{cor}[prop]{Corollary}
\newtheorem{lem}[prop]{Lemma}
\newtheorem{problem}[prop]{Problem}
\newtheorem{question}[prop]{Question}
\newtheorem{con}[prop]{Conjecture}

\theoremstyle{definition}
\newtheorem{defn}[prop]{Definition}
\newtheorem{ack}{Acknowledgments} 
\theoremstyle{remark}
\newtheorem{rem}{Remark}

\newcommand{\bbN}{{\Bbb{N}}}

\newcommand{\bbZ}{{\Bbb{Z}}}

\newcommand{\calC}{{\cal{C}}}

\newcommand{\calE}{{\cal{E}}}

\newcommand{\calP}{{\cal{P}}}
\newcommand{\calR}{{\cal{R}}}
\newcommand{\al}{\alpha}
\newcommand{\be}{\beta}
\newcommand{\g}{\gamma}

\newcommand{\e}{\varepsilon}

\newcommand{\la}{\lambda}

\renewcommand{\span}{\operatorname{span}}
\newcommand{\supp}{\operatorname{supp}}
\newcommand{\card}{\operatorname{card}}
\newcommand{\conv}{\operatorname{conv}}
\newcommand{\codim}{\operatorname{codim}}

\newcommand{\Ker}{\operatorname{Ker}}
\newcommand{\sgn}{\operatorname{sgn}}
\renewcommand{\Re}{\operatorname{Re}}

\newcommand{\disp}{\displaystyle}
\newcommand{\lb}{\label}

\newcommand{\lra}{\longrightarrow}

\newcommand{\wtw}{if and only if }

\newcommand{\ONTO}{\buildrel {\mbox{\small onto}}\over \longrightarrow}
\newcommand{\INTO}{\buildrel {\mbox{\small into}}\over \longrightarrow}

\makeatletter
\def\@currentlabel{2.1}\label{e:dispaa}
\def\@currentlabel{2.21}\label{e:dispau}
\def\@currentlabel{2.22}\label{e:dispav}
 \def\@currentlabel{2.23}\label{e:dispaw}
\def\@currentlabel{2.24}\label{e:dispax}
\makeatother

\makeatletter
\def\alphenumi{%
  \def\theenumi{\alph{enumi}}%
  \def\p@enumi{\theenumi}%
  \def\labelenumi{(\@alph\c@enumi)}}
\makeatother

\begin{document}

\title{Norm One Projections in Banach Spaces}

\author{Beata Randrianantoanina}

\address{Department of Mathematics and Statistics
\\ Miami University \\Oxford, OH 45056, USA}

 \email{randrib@@muohio.edu}

\begin{abstract}
This is the survey of results about norm one projections and 
$1$-complemented subspaces in K\"othe function spaces and Banach 
sequence spaces.  The historical development of the theory is 
presented from the 1930's to the newest ideas. Proofs of the main 
results are outlined.  Open problems are also discussed.  Every 
effort has been made to include as complete a bibliography as 
possible.
\end{abstract}


\maketitle

\tableofcontents

\section{Introduction}

One of the main topics in the study of Banach spaces has been, since the
inception of the field, the study of  projections and complemented
subspaces. Here by a projection we mean a bounded linear operator $P$
satisfying $P^{2}=P$, and by a complemented subspace we mean a range of a
bounded linear projection $P$.

Banach \cite{Ban} posed several problems about  projections and complemented
subspaces.  Some of the most famous ones are:

\begin{problem} \lb{basis}
Does every complemented subspace of a 
space with a basis 
have a basis?
\end{problem} 

\begin{problem}  \lb{comp} Does every complemented subspace of a space with an 
unconditional
basis have  an unconditional basis?
\end{problem} 

\begin{problem} \lb{bpro} Does every Banach space $X$ admit a nontrivial bounded linear
projection? 
\end{problem} 

Here nontrivial means that $\dim P(X)= \infty$ and
$\dim X/P(X)=\infty$; notice that the Hahn-Banach theorem guarantees that for
every Banach space $X$ and every subspace $Y \subset X$ with 
$\dim Y < \infty$
there exists a bounded linear projection with $P(X)=Y$.

These problems proved to be very elusive.  In fact Problem~\ref{comp} is still
open, even for subspaces of $L_{p}$, despite the intensive work of many
people in the area.

Problems~\ref{basis} and \ref{bpro} were both answered negatively more 
than 50 years after
they were posed.  Problem~\ref{basis} was solved by Szarek in 1987 
\cite{Szarek87}
and Problem~\ref{bpro} was solved by Gowers and Maurey in 1993 \cite{GM93}.

The work, of more than 50 years, on Problems~\ref{basis}, \ref{comp}, 
\ref{bpro} has 
 led to significant
developments in Banach space theory and also to many intriguing
questions.

In this survey we want to concentrate on the developments of the theory of
projections of norm one in Banach spaces.  It seems that this theory, being
of isometric rather than isomorphic nature, should be much less
complicated than the theory of bounded projections of arbitrary norm.

Note, however, that any Banach space can be equivalently 
renormed so that the given
bounded projection $P$ has norm one on $X$.  Thus, without loss of
generality, one can rephrase Problem~\ref{basis} and \ref{comp} to ask for
the $1$-complemented subspaces with the additional properties as stated
there.  Maybe this is the reason why the theory of norm one projections
still has many open problems.

One of the most interesting of them is the isometric version of
Problem~\ref{comp}:

\begin{problem}
Does every $1$-complemented subspace of a space $X$ with a $1$-unconditional
basis have  an unconditional basis (with any constant $C$)?
\end{problem}

This problem has the affirmative answer if space $X$ is over $\Bbb C$
\cite{KW} (and then constant $C=1$), but it is open if space $X$ is real, see
Section~\ref{sequence} (it is known that $C=1$ does not
 work in the real case).
We will discuss also some other interesting open problems related to
Problem~\ref{bpro}, see Questions~\ref{renorm}, \ref{subspace}.

Norm one projections are important because they are one of the most natural
generalizations of the concept of orthogonal projections from Hilbert
spaces to arbitrary Banach spaces.  

Another very natural generalization of orthogonal projections are metric projections
 (also called proximity
mappings, nearest point mappings, cf.
Definition~\ref{defprox}), which play a very important role in the
theory of best approximation.

Metric projections are usually set-valued and one of the  major research
directions in the area is to determine when does a metric projection admit
a continuous or linear selection (see e.g. \cite{S74,Deutsch83}).  Not
surprisingly contractive projections and metric projections as two natural
generalization of orthogonal projections are
intrinsically related to each other  and there
are results about linear selections of metric projections using
characterizations of norm-one projections as well as results giving
characterizations of contractive projections using facts about metric
projections (see Section~\ref{metric}).

It is striking that despite the intensive work of many authors on
contractive projections (our bibliography includes over 120 items and
probably is not complete), the full characterization of norm one projections
 is known only in
$L_p$-spaces among nonatomic K\"othe function spaces.

There are examples demonstrating that $1$-complemented subspaces of
general Banach spaces other than $L_p$ cannot be described purely in
isomorphic or isometric terms (in $L_p$ a subspace $X$ is $1$-complemented
if and only if $X$ is isometrically isomorphic to an $L_p$-space, possibly
on a different measure space).  Thus it seems that the $1$-complementability
of a subspace $Y$ in $X$ in fact depends on the way that $Y$ is embedded in
$X$ and that norm one projections are best described in terms of 
conditional expectation operators in both nonatomic K\"othe function spaces
and sequence Banach spaces, see
 Conjectures~\ref{condexpcon} and \ref{blockcon}.

This survey is organized as follows.  We start from recalling  
well-known facts about the abundance of contractive projections in Hilbert
spaces and the characterizations of Hilbert spaces through existence of
enough contractive projections (Section~\ref{Hilbert}).

Next we divide the survey into two parts.  Part I is about nonatomic function
spaces and Part II is about sequence spaces.  Each part starts with the
section devoted to Lebesgue spaces $L_p$ and $ \ell_p$, respectively
(Sections~\ref{Lpsection}, \ref{lpsection}).  
  Sections~\ref{function} and \ref{sequence} are devoted to a variety of
partial results valid in K\"othe function spaces and in sequence spaces.

In the present survey we limit ourselves to the theory of contractive
projections in K\"othe function spaces and sequence spaces, which are not
$M$-spaces.  There exists vast literature on contractive projections in
spaces of continuous functions, vector valued function spaces,
noncommutative Banach spaces as well as nonlinear contractive projections,
so we could not possibly cover everything here.  We refer the interested
reader to surveys \cite{PelB,Doust95} of some of these topics, some
additional references are also mentioned in the text below.

Throughout the survey we use standard notations as may be found e.g. in
\cite{LT1,LT2}.  For the convenience of the reader we collect in 
Section~\ref{Appendix} some of the most important definitions that we use.

\begin{ack}
I want to thank the organizers of the ICMAA2000 for giving me the
opportunity to speak at this nice conference and for their hospitality
during my stay in Kaoshiung. I thank Miami University for the funding 
for my travel to this conference.

I also want to thank Ms. L. Ferriell for the patient typing of this
manuscript.
\end{ack}

\section{Definitions}\lb{Appendix}

As indicated in the Introduction throughout the paper we use standard
definitions and notations as may be found e.g. in \cite{LT1,LT2}.  In this
section, for the convenience of the reader, we collect some of the most
important definitions that we use. Throughout the main body of the paper
we will refer to these definitions whenever we use them. 

\begin{defn}
\cite[Definition~1.b.17]{LT2}
\lb{defkfs}
Let $(\Omega,\Sigma, \mu)$ be a complete $\sigma$-finite measure space.  A
Banach space $X$ consisting of equivalence classes, modulo equality almost
everywhere, of locally integrable functions on $\Omega$ is called a {\it
K\"othe function space} if the following conditions hold:
\begin{itemize}
\item[(1)]
If $|f(w)| \leq |g(w)|$ a.e. on $\Omega$, with $f$ measurable and $g \in
X$, then $f \in X$ and $\|f \|_X \leq \|g\|_X$,
\item[(2)]
For every $A \in \Sigma$ with $\mu(A)<\infty$ the characteristic function
$\chi_A$ of $A$ belongs to $X$.
\end{itemize}
\end{defn}

\begin{defn}
\lb{defoc}
We say that a K\"othe function space $X$ is {\it order-continuous} if
whenever $f_n \in X$ with $f_n \downarrow 0$ a.e. then $\| f_n \|_X
\downarrow0$.
\end{defn}

\begin{defn}
\lb{defkdual}
The K\"othe dual $X'$ of $X$ is the K\"othe space of all $g$ such that
$\int_\Omega |f| \ |g| \ d \mu < \infty$ for every $f \in X$.  $X'$ is
equipped with the norm
$$
\| g \|_{X'} = \sup \limits_{\|f\|_{X}\leq 1} \ \int \limits_\Omega |f| \ |g|
d \mu.
$$
\end{defn}
K\"othe dual $X'$ can be regarded as a closed subspace of the dual $X^*$ of
$X$. If $X$ is order-continuous then $X' = X^*$.

\begin{defn}
\cite[Definition~2.a.1]{LT2} \lb{defri}
Let $(\Omega, \Sigma, \mu)$ be one of the measure spaces
$[0,1]$, $[0,\infty)$ or $\{1,2,\dots\}$ (with the natural measure).  A
K\"othe function space $X$ on $(\Omega, \Sigma, \mu)$ is said to be a {\it
rearrangement invariant space} ({\it r.i. space}) if the following
conditions hold:

\begin{itemize}
\item[(1)]
If $\tau$ is a measure-preserving automorphism of $\Omega$ onto itself and
$f$ is measurable function on $\Omega$ then $f \in X$ if and only if $f
\circ \tau \in X$ and $\|f \|_X = \|f \circ \tau \|_X$.
\item[(2)]
$X'$ is a norming subspace of $X^*$.
\item[(3)]
If $A \in \Sigma$ and $\mu(A)=1$ then $\| \chi_A \|_X = 1$.
\end{itemize}
\end{defn}

R.i. spaces are also sometimes called {\it symmetric}, especially when
 $\Omega=\{1,2,\dots\}$.  

The most commonly
used r.i. spaces, besides the Lebesgue spaces $L_p$ and $\ell_p$, are 
Orlicz   and Lorentz spaces.

\begin{defn} \lb{defo}
(see e.g. \cite[Definition~4.a.1]{LT1})
An {\it Orlicz  function} $\varphi$ is a left-con\-ti\-nu\-ous, non-decreasing
  convex function $\varphi:[0, \infty) \lra [0,\infty]$ such that
$\varphi(0) =0$ and $\lim_{t \to \infty} \varphi (t) = \infty$.  We
will also additionally assume that $\varphi(1) = 1$.

Let $(\Omega, \Sigma, \mu)$ be one of the measure spaces $\{1,2,
\dots\}, \ [0,1]$ or $[0, \infty)$ (with the natural measure).  The Orlicz
spaces $L_\varphi(\Omega, \Sigma, \mu)$ is the space of all (equivalence
class of) measurable functions $f$ on $\Omega$ so that
$$
\int \limits_\Omega \varphi (\frac {|f(t)|} {\lambda}) d \mu < \infty
$$
for some $\lambda>0$.  The space $L_\varphi$ is equipped with the norm
$$
\|f \|_\varphi = \inf \{\lambda > 0: \int \limits_\Omega
\varphi(\frac{|f(t)|}{\lambda}) \leq 1\}.
$$

$\| \cdot \|_\varphi$ is called a {\it Luxemburg norm}.

Sometimes $L_{\varphi}$ is considered with a different equivalent norm
$\Vert \cdot \Vert_{\varphi, O}$ which can be described as follows:

\begin{equation}\label{Amemiya}
\Vert x \Vert_{\varphi,O} = \inf \limits_{\lambda > 0} \frac{1}{\lambda}
(1 + \int \limits_\Omega\varphi 
(\lambda \vert x_{n} \vert)d\mu)
\end{equation}

$\Vert \cdot \Vert_{\varphi,O}$ is called an {\it  Orlicz  norm} and
$\eqref{Amemiya}$ is called an {\it Amemiya  formula  
for  the Orlicz 
norm}.
\end{defn}

Orlicz functions were introduced   by Birnbaum and Orlicz
\cite{BO31} and Orlicz spaces were first considered by Orlicz
\cite{Or32,Or36}. Since then they were extensively studied by many authors
and became a source of many examples and counterexamples in Banach space
theory.  There are many monographs devoted to Orlicz spaces (see e.g.
\cite{Kr-Rut,RaoRen,Chen}).

The other natural extension of Lebesgue spaces $L_{p}$ and $\ell_p$ 
are Lorentz
 spaces $L_{w,p}$ and  $\ell_{w,p}$ introduced by Lorentz \cite{Lor50} in
connection with some problems of harmonic analysis and interpolation
theory. In this paper we do not use Lorentz function spaces, so we just 
recall the
definition of Lorentz sequence spaces $\ell_{w,p}$.

\begin{defn}
(see e.g. \cite[Definition~4.e.1]{LT1})
\lb{defl}
Let $1 \leq p < \infty$ and let $w = \{w_{n}\}^{\infty}_{n=1}$
be a non-increasing sequence of non-negative numbers such that
$w_1=1$ and $\lim_{n \to \infty} w_{n}=0$.  The
Banach space of all sequences of scalars 
$x=(x_{1}, x_{2},\dots)$ for
which
$$
\Vert x \Vert_{w,p} = \sup_{\pi}
(\sum^{\infty}_{n=1} \vert a_{\pi (n)} \vert^{p}
w_{n})^{\frac{1}{p}} < \infty
$$
where $\pi$ ranges over all permutations of the integers, 
is called a {\it
Lorentz  sequence  space} and is denoted by $\ell_{w,p}$
(or $d(w,p)$).
\end{defn}

Note that $\Vert \cdot \Vert_{w,p}$ can also be computed as follows:
$$
\Vert x \Vert_{w,p} =(\sum_{n=1}^{\infty}(x_{n}^{*})^{p}
w_{n})^{\frac {1}{p}}
$$
where $x=(x_{1}, x_{2} \cdots) \in \ell_{w,p}$ and
$\{x_{n}^{*}\}^{\infty}_{n=1}$ is a non-increasing sequence obtained from
$\{\vert x_{n} \vert \}^{\infty}_{n=1}$ by the suitable permutation of the
integers.  $\{x_{n}^{*}\}^{\infty}_{n=1}$ is called a 
{\it non-increasing
 rearrangement  of  the  sequence}  $\{x_{n}\}^{\infty}_{n=1}$.

\begin{defn}
\lb{defJ}
Let $X$ be a Banach space.  We define a {\it duality map} $J$ from $X$ into
subsets of $X^*$ by the condition that $f \in J(x)\subset X^*$ if and only
if $ \|f \|_{X^{*}} = \| X \|_X$ and $\left< f, x \right> = \|x \|^2_X$.

If $J(x)$ contains exactly one functional then element $x$ is called {\it
smooth} in $X$.

If every element $x\in X$ is smooth in $X$ then $X$ is called {\it
smooth}.
\end{defn}

\begin{defn}
\lb{defsm}
We say that the norm in the Banach lattice $X$ is {\it strictly monotone}
if $\|x +y\|_X > \|x \|_X$  whenever $x, y \geq 0$ and $y \neq 0$.
\end{defn}

\begin{defn} \lb{defbasis}
A Schauder basis $\{x_i\}_i$ for $X$ is called {\it monotone} if $\sup
\limits_n \|P_n \| =1$ where $P_n$ are the natural projections associated
to the basis i.e.
$$
P_n(\sum \limits^\infty_{i=1} a_i x_i)= \sum \limits^n_{i=1}a_i x_i.
$$

If $\sup \limits_n \|I-P_n\| =1$  (where $I$ denotes the identity operator)
the basis is called {\it reverse monotone}.

 If for all scalars $(a_i)_i$ we
have
$$
\| \sum \limits^{\infty}_{i=1} a_i x_i \|_X = \| \sum \limits^\infty_{i=1}
| a_i | x_i \|_X
$$
then the basis $\{x_i\}_i$ is called {\it $1$-unconditional}.
\end{defn}

Most commonly studied  examples of spaces with
$1$-unconditional bases include $\ell_p$, Orlicz and Lorentz spaces.

\begin{defn}
\lb{defprox}
Let $X$ be a real or complex normed linear space and $M$ a subset of $X$.
The {\it metric projection} onto $M$ is the mapping $P_M: X \to 2^M$ which
associates with each $x$ in $X$ its (possibly empty) set of nearest points
in $M$ i.e.
$$
P_M(x)=\{m \in M: \|x-m\|= \inf \{\|x-y \|: y \in M\}\}
$$

Other terms for metric projections used in the literature include {\it best
approximation operator, nearest point map, Chebyshev map, proximity
mapping, normal projection, projection of minimal distance}.
\end{defn}

Some authors use the term metric projection for a particular selection of
the set-valued mapping $P_M$ (cf. Definition~\ref{defsel}).

\begin{defn}
\lb{defCh}
A subset $M$ of a normed linear space $X$ is called {\it proximinal}
(resp. {\it Chebyshev}) if for each $x \in X, \ P_M(x)$ contains at least one
(resp. exactly one) element of $M$.
\end{defn}

\begin{defn}
\lb{defsel}
Let $Q: X \to 2^Y$  be a set-valued map.  A {\it selection} for $Q$ is any
map $q: X \to Y$ such that $q(x)\in Q(x)$ for all $x \in X$.
\end{defn}

\section{Hilbert spaces} \lb{Hilbert}
The first results about norm one projections were proven in the setting of
Hilbert spaces.

It is well-known that Hilbert spaces contain many norm one projections.
Namely, orthogonal projections are contractive and for every subspace
$Y$ of a Hilbert space $H$ there exist an orthogonal projection whose range
is precisely $Y$.  In fact this property characterizes Hilbert spaces as
was proven by Kakutani \cite{Ka39} (see also \cite{Ph40}) in the case of
real spaces, and by Bohnenblust \cite{Bh42} in the complex case (cf.
\cite{PelB}). This was later refined by James \cite{J47} and Papini 
\cite{Pap74} (see also \cite{Reich83}). We have: 

\begin{th} \label{Hspace}
(cf. \cite {Amir}) For a Banach space $X$ with $\dim X \geq 3$
the following statements are equivalent:

\begin{itemize}
\item[(i)]
$X$ is isometrically isomorphic to a Hilbert space,

\item[(ii)]
every 2-dimensional subspace of $X$ is the range of a projection of norm 1,

\item[(iii)]
every subspace of $X$ is the range of a projection of norm 1,

\item[(iv)]
(James \cite{J47}) every 1-codimensional subspace of $X$ is the range of a
projection of norm $1$,

\item[(v)]
(Papini \cite{Pap74}, de~Figueiredo, Karlovitz \cite{FK} for the case when
$\dim X<\infty$ and $X$ is strictly convex) for some $2 \leq n \leq \dim X-1$, every
n-dimensional subspace of $X$ is the range of a projection of norm 1,

\item[(vi)]
(Papini \cite{Pap74}) for some $1 \leq n \leq \dim X-2$, every
n-codimensional subspace of $X$ is the range of a projection of norm $1$.
\end{itemize}
\end{th}

Note also that in Hilbert spaces every subspace is isometrically isomorphic
to a Hilbert space and therefore the equivalence $(i)\leftrightarrow (iii)$
in Theorem~\ref{Hspace} can be restated as follows:

\begin{prop} \lb{isometryH}
Let $X$ be a Hilbert space.  Then $Y \subset X$ is contractively
complemented if and only if $Y$ is isometrically isomorphic to a Hilbert
space.\end{prop}

However, as we will see in the sequel, the statement in  
Proposition~\ref{isometryH} does not
characterize Hilbert spaces.
\part{Nonatomic function spaces}

\section{Lebesgue function spaces ${L_p}$}\lb{Lpsection}

It follows from Theorem~\ref{Hspace} that in spaces other than Hilbert
space there are many subspaces which are not $1$-complemented.  But even
before the results of Kakutani (1939) and Bohnenblust (1942), Murray
\cite{M37} showed, answering a question of Banach \cite{Ban}, that if
$1<p<\infty,\  p \neq 2$ then there exist subspaces of $L_{p}$ and $\ell_{p}$
which are not complemented.  This appears to be the first result indicating
that not every subspace of $L_{p}$ is $1$-complemented.

\subsection{Conditional expectation operators}

In 1933, in the treatise on the foundations of probability 
\cite{Kol33}, 
Kolmogorov introduced  conditional expectation operators, which are very 
important examples of norm one projections in $L_p$ and other Banach spaces.
Moreover
conditional expectation operators play crucial role in describing general
contractive projections, so we  start from recalling
their definitions and basic properties.

\begin{defn} \cite{Kol33,Doob53}
Let $(\Omega, \Sigma, \mu)$  be a $\sigma$-finite measure space and
$\Sigma_{0} \subset \Sigma_{0}$ be a $\sigma-$subalgebra of $\Sigma$
such that $\mu$ restricted to $\Sigma_{0}$ is $\sigma-$finite (i.e. so
that $\Sigma_{0}$ does not have atoms of infinite measure).  By the
Radon-Nikodym Theorem, for every
$
f \in L_1 (\Omega, \Sigma, \mu) +L_{\infty}(\Omega, \Sigma, \mu)
$
there exists a unique, up to an equality a.e., $\Sigma_{0}-$measurable
locally integrable function $\calE^{\Sigma_{0}}f$ so that
$$
\int_{\Omega} g \calE^{\Sigma_{0}} f d \mu = \int_{\Omega} g f d \mu
$$
for every bounded, integrable, $\Sigma_{0}-$measurable function $g$.
The function $\calE^{\Sigma_{0}} f$ is called the {\it conditional 
expectation  of  $f$  with  respect  to  $\Sigma_{0}$}.
\end{defn}

Basic properties of conditional expectations are as follows
\cite{Kol33,Doob53}:
\begin{itemize}
\item[(1)]
$\calE^{\Sigma_{0}}$ is a linear mapping,
\item[(2)]
$\calE^{\Sigma_{0}}$ is positive i.e. $\calE^{\Sigma_{0}} (f) \geq 0$
whenever $0 \leq f \in L_1(\Omega, \Sigma, \mu) + L_{\infty}
(\Omega, \Sigma, \mu)$,
\item[(3)]
$\calE^{\Sigma_{0}} (f)=f$  if and only if
$
f \in L_1(\Omega, \Sigma_{0}, \mu) + L_{\infty} (\Omega, \Sigma_{0},
\mu),
$
in particular $\calE^{\Sigma_{0}} (\bold1)  =  \bold1$ (here
$\bold1$ denotes a
function on $\Omega$ constantly equal to $1$),
\item[(4)]
$\Vert \calE^{\Sigma_{0}} f \Vert_{1} \leq \Vert f \Vert_{1}$  for all
$f \in L_1 (\Omega, \Sigma, \mu)$ (this follows easily since if $f
\geq 0$ then $\Vert \calE^{\Sigma_{0}} f \Vert_{1} = \Vert f \Vert_{1}$
and by positivity of $\calE^{\Sigma_{0}}$ we have $\vert
\calE^{\Sigma_{0}} (f) \vert \leq \calE^{\Sigma_{0}}(\vert f \vert))$
which gives the desired conclusion),
\item[(5)]
(cf. \cite{Moy54})
if $p \geq 1$ then $\vert \calE^{\Sigma_{0}} (f) \vert^{p} \leq
\calE^{\Sigma_{0}} (\vert f \vert^{p})$ almost everywhere so
$ \Vert \calE^{\Sigma_{0}} (f) \Vert_{p} \leq \Vert f \Vert_{p}$
for all $f \in L_{p} (\Omega, \Sigma, \mu)$,
\item[(6)]
$\calE^{\Sigma_{0}}$ satisfies the ${\it averaging \ identity}$ i.e.
$$
\calE^{\Sigma_{0}} (f \cdot \calE^{\Sigma_{0}} (g)) =
\calE^{\Sigma_{0}}(f) \cdot \calE^{\Sigma_{0}} (g)
$$
for all $f \in L_{\infty} (\Omega,\Sigma,\mu)$ and $g \in
L_1(\Omega, \Sigma, \mu) + L_{\infty} (\Omega, \Sigma, \mu)$.
\end{itemize}

These properties can be summarized as follows:
\begin{th}  \lb{condexp}
(for $p=1$ implicit in \cite{Kol33}, for $p>1$ see \cite{Doob53,Moy54})
If $p \geq 1$ and $(\Omega, \Sigma, \mu)$ is a $\sigma$-finite measure space and $\Sigma_{0}$
is a $\sigma$-subalgebra of $\Sigma$ then the conditional expectation
operator $\calE^{\Sigma_{0}}$ is a positive contractive projection from
$L_{p}(\Omega, \Sigma, \mu)$ onto $L_{p} (\Omega, \Sigma_{0}, \mu)$, that
leaves the constants invariant.
\end{th}

Several authors have investigated whether any subset of
properties $(1)-(6)$ characterizes conditional expectation operators.  For
the most recent results in this direction as well as a nice account of the
literature we refer the reader to \cite{DHP90}.  

In this survey we will
concentrate on the papers that deal directly with the converse of Theorem~\ref{condexp}.

First such results go back to the 50's
\cite{Moy54,Bahadur55,Sidak57,Rota60,Brunk63}.  These papers dealt with
converse of Theorem~\ref{condexp} under some additional assumptions on the
contractive projection operator.

\subsection{Main characterizations}
In 1965 Douglas \cite{D}  gave the complete
characterization of contractive projections on $L_1 (\Omega,
\Sigma,\mu)$ when $(\Sigma,\mu)$ is a finite measure space.  In the following year Ando
\cite{Ando} extended this characterization to $L_{p}(\Omega,\Sigma,\mu)$ on
a finite measure space for arbitrary $p, 1 \leq p < \infty,\  p \neq 2$.
They proved that every norm one projection on $L_{p}(\Omega, \Sigma, \mu)$,
$1 < p < \infty, p \neq 2$, which leaves constants invariant is a
conditional expectation operator (when $p=1$ this is true for most norm
one projections, for the precize statement see the
 following Theorem~\ref{ando}). Moreover if the norm one projection $P$
does not leave constants invariant, i.e. if $P{\bold{1}}= h$ is an arbitrary function
in $L_{p}(\Omega,\Sigma,\mu)$ then $h$ has a maximum support among functions 
from the range of $P$, the range of $P$ is isometrically isomorphic to a 
weighted $L_p-$space (with weight $h$) and the projection $P$ is a, so called,
weighted conditional expectation operator. The precize statement of their
results is given below:

\begin{th} \cite{D,Ando}\lb{ando}
Let $1 \leq p < \infty, p \neq 2$ and let $(\Omega, \Sigma, \mu)$ be a
finite measure space.  Then $P: L_{p}(\Omega, \Sigma, \mu) \longrightarrow
L_{p}(\Omega, \Sigma, \mu)$ is a contractive projection if and only if
there exists a $\sigma$-subalgebra $\Sigma_{0} \subset \Sigma$ and a
function $h \in L_{p}(\Omega, \Sigma, \mu)$ such that the support $B$ of
$h$ is the maximum element of $\Sigma_{0}$ and 
 for all $f \in L_{p}(\Omega, \Sigma, \mu)$, $P$ is represented in the
form:
$$
Pf = \frac{h}{\calE^{\Sigma_{0}} (\vert h \vert^{p})} \cdot
\calE^{\Sigma_{0}} (f \cdot \bar{h}^{p-1}) + Vf
$$
where, when $p \neq 1, V = 0$
and when $p=1$, $V$ is a contraction such that $V^{2}=0, P_{B}V=V,
VP_{B}=0$ and $ {Vf}/{h}$ is $\Sigma_{0}-$measurable (here $P_{B}$
denotes the projection defined by $P_{B}f = f \cdot \chi_{B}$, where
$\chi_{B}$ is the characteristic function of set $B$).
\end{th}

The method of Douglas and Ando depends on studying positive projections
and projections whose range is a sublattice.  Subsequently, as we will
discuss in Section~\ref{function}, many of the key steps leading to
Theorem~\ref{ando} were generalized to other spaces.  Thus we will
present these key steps here:

\begin{prop} \lb{latrep}
(related to results of \cite{Bahadur55,Brunk63,Moy54})
Suppose that $Y \subset L_{1}$ is a closed sublattice of $L_{1}(\Omega,
\Sigma, \mu)$.  Then there exists a $\sigma$-subalgebra $\Sigma_{0} \subset
\Sigma$ and a weight function $k$ such that $Y=k \cdot L_{1} (\Omega,
\Sigma_{0}, \mu_{\Sigma_{0}})$.  Moreover the pair $(\Sigma_{0},k)$ is
unique.  (The measure $\mu_{\Sigma_{0}}$ is the restriction of $\mu$ to
$\Sigma_{0}$.)
\end{prop}

\begin{prop} \lb{lattice}
If $P$ is a positive contractive projection on $L_{1}(\Omega,\Sigma, \mu)$,
then $\calR(P)$ is a closed sublattice of $L_{1}(\Omega, \Sigma, \mu)$.
(Here $\calR(P)$ denotes the range of $P$.)
\end{prop}

This follows from the fact that if $f \in \calR(P)$ then we have
$$
f^{+} \geq f \Longrightarrow P(f^{+}) \geq P(f) =f
$$
and therefore $P(f^{+}) \geq f^{+} \geq 0$.

Thus
$$
\Vert P(f^{+}) - f^{+} \Vert_{1} = \int_{\Omega} (P(f^{+})-f^{+}) d \mu =
\Vert P(f^{+}) \Vert_{1} - \Vert f^{+} \Vert_{1} \leq 0.
$$

Hence $\Vert P(f^{+})- f^{+} \Vert_{1} = 0$ i.e. $P(f^{+}) = f^{+}$ so
$f \in \calR(P) \Longrightarrow f^{+} \in \calR(P)$ and the conclusion of
Proposition~\ref{lattice} follows easily.

\begin{th} \lb{P1=1}

For a linear operator $P: L_{1}(\Omega, \Sigma, \mu) \to L_{1}(\Omega,
\Sigma, \mu)$  the following are equivalent:

\begin{itemize}
\item[(a)] $P$ is a contractive projection with $P (\bold 1) = \bold 1$

\item[(b)] $P$ is a conditional expectation operator, i.e. there exists a
unique $\sigma$-subalgebra $\Sigma_{0} \subset \Sigma$ such that
 $P = {\cal E}^{\Sigma_0}$
\end{itemize}

\end{th}
Very elegant and short proof of this theorem was given recently by
Abramovich, Aliprantis and Burkinshaw \cite{AAB93}, so we will not reproduce
it here.  
The main idea is to
 first show  that every contraction $T$ in
$L_{1}$ with $T(\bold 1) = \bold 1$ is positive (first proven by Ando
\cite{Ando}) and then use Propositions~\ref{lattice} and \ref{latrep} to
see that $\calR(P)=L_{1}(\Omega, \Sigma_{0}, \mu_{\Sigma_{0}})$ for some
$\sigma$-subalgebra $\Sigma_{0} \subset \Sigma$ (with $\Omega \in
\Sigma_{0}$).  In particular we see that $P\chi_{A} =\chi_{A}$ for all $A \in
\Sigma_{0}$, and Theorem~\ref{P1=1}  quickly follows (see \cite{AAB93}).

Next Douglas notices that in $L_{1}$ there exist contractive projections of
a certain ``irregular'' form, which is impossible in $L_{p}$ with $p>1$;
namely we have:

\begin{lem} \lb{L1V}
Let $\Sigma_{0} \subset \Sigma$ be a $\sigma$-subalgebra with a maximal
element $B  \subsetneq \Omega$.  Let $V:L_{1} (\Omega, \Sigma, \mu) \to
L_{1}(\Omega,\Sigma,\mu)$ be a contractive linear operator such that
$V^{2}= 0,  P_{B}V = V$ and $VP_{B}=0$.
Then the operator $P={\cal E}^{\Sigma_0}+V$ is a contractive linear
projection from $L_{1} (\Omega, \Sigma, \mu)$ onto $L_{1} (B,
\Sigma_{0},\mu_{\Sigma_{0}})$.
\end{lem}

\begin{pf*}{Sketch of Proof}
It is not difficult to check that $P^{2}=P$ and
\begin{align*}
Pf &= P(f\chi_{B} + f\chi_{B^{c}}) \\
&= {\cal E}^{\Sigma_0} (f\chi_{B}) + V(f\chi_{B}) + {\cal E}^{\Sigma_0} (f\chi_{B^{c}})
+ V(f\chi_{B^{c}}) \\
&= {\cal E}^{\Sigma_0} (f\chi_{B}) + 0 + 0 + V (f\chi_{B^{c}}).
\end{align*}

Thus
\begin{align*}
\Vert P{f} \Vert_{1} &\leq \Vert {\cal E}^{\Sigma_0} (f\chi_{B}) \Vert_{1} +
\Vert V (f\chi_{B^{c}}) \Vert_{1} \\
&\leq \Vert f\chi_{B} \Vert_{1} + \Vert f\chi_{B^{c}} \Vert_{1} = 
\Vert f \Vert_{1}\ .
\end{align*}

So $P$ is contractive.
\end{pf*}

Moreover, Ando showed:
\begin{th} \lb{reductiontoL1}
A contractive projection $P$ in $L_{p} \ (1 < p < \infty,\  p \neq 2)$ with $P
\bold 1 = \bold 1$ is contractive with respect to $L_{1}$-norm.
\end{th}

The final step in the proof of Theorem~\ref{ando} is to drop the assumption
that $P \bold 1 = \bold 1 $ in Theorem~\ref{P1=1}.
(So far any attempts to generalize this final step to spaces other than
$L_{p}$, have failed, see the discussion in Section~\ref{function}).

Theorem~\ref{ando} can be restated as follows:

\begin{th} \lb{isometryequiv}
\cite[Corollary~1]{Ando}
A contractive projection $P$ in $L_{p}$ is isometrically equivalent to a
conditional expectation (with respect to a measure), if $1<p<\infty$, or if
$p=1$ and $PP_{B}=P$, where $B$ is the maximum support of the range of
$P$.
\end{th}

Also as a corollary one immediately obtains an analogue of
Proposition~\ref{isometryH}.  Namely we have:

\begin{th}
\lb{isometryLp}
(Douglas, Ando, see \cite[Theorem~4.1]{BL74} for general measures)

\begin{itemize}
\item[(A)]
For a subspace $Y$ of $L_{p} (\Omega, \Sigma, \mu)$ the following
statements are equivalent:
\begin{itemize}
\item[(A1)]
$Y$ is the range of a positive contractive projection,

\item[(A2)]
$Y$ is the closed sublattice of $L_{p}(\Omega, \Sigma, \mu)$,

\item[(A3)]
there exists a positive isometrical isomorphism from $Y$ onto some
$L_{p}(B, \Sigma_{0}, \nu)$.
\end{itemize}
\item [(B)]
For a subspace $Y$ of $L_{p}(\Omega, \Sigma, \mu)$ the following statements
are equivalent:
\begin{itemize}
\item[(B1)]
$Y$ is the range of a contractive projection,

\item[(B2)]
there exists function $\phi$ on $\Omega$, with $\vert \phi \vert =1$ a.e.
such that $\phi \cdot Y$ is the closed sublattice of $L_{p} (\Omega,
\Sigma, \mu)$,

\item[(B3)]
$Y$ is isometrically isomorphic to some $L_{p} (B, \Sigma_{0}, \nu)$.
\end{itemize}
\end{itemize}
\end{th}

\subsection{Characterizations of $L_p$ through contractive projections}
In 1969 Ando \cite{ando69} showed that Theorem~\ref{isometryLp} does
characterize $L_{p}$ among Banach lattices of dimension bigger or equal
than $3$.  He proved:

\begin{th}
\lb{andolattice}

Let $X$ be a Banach lattice with $ \dim X \geq 3$.  Then $X$ is order
isometric to $L_{p}(\mu)$, for some $1 \leq p < \infty$ and measure $\mu$,
or to $c_{0}(\Gamma)$, for some index set $\Gamma$, if and only if there is
a contractive positive projection from $X$ onto any closed sublattice of
it.
\end{th}

For expositions of this theorem see e.g.
\cite{BLcharLp}, \cite[\S 16]{Lacey}, \cite[Section~1.b]{LT2}.
Theorem~\ref{andolattice} has been strengthened in a series of articles by 
Calvert and Fitzpatrick (1986--1991) which will be discussed in 
Section~\ref{lpsection}.
Here we just mention one of their theorems which directly applies to
nonatomic function spaces.  

As Calvert and Fitzpatrick indicate, in the
following statement, they had in mind $A$ being a set of characteristic
functions of measurable sets:

\begin{th}
\cite[Theorem~4.4]{CF87} \lb{CFnonatomic}
Let $X$ be a Banach lattice and suppose $A = \{e_i\}_{i \in I}$ is a set of
positive elements of $X$ such that $e,f \in A$ implies $(e-f)^+ \in A$ and
$(e-f)^+ \wedge f = 0$.  Suppose $X$ is the closed span of $A$, and $\dim
(\span A) \geq 3$.  Suppose that any two dimensional sublattice of $X$
containing any $e_i$, $i \in I \backslash  \{i_0 \}$ (where $i_0 \in I$ is a fixed
arbitrary element of $I$) is the range of a contractive projection.  Then
$X$ is an $L_p$-space $(1 \leq p< \infty)$ or an $M$-space.
\end{th}

\subsection{Extensions of the main characterizations to general measure 
spaces}
Earlier, in 1955, Grothendieck \cite{G55}  
showed part of Theorem~\ref{isometryLp}, namely he showed 
that when $p=1$ and
$(\Omega, \Sigma, \mu)$ is a general measure space then a range of a
contractive projection in $L_1(\Omega, \Sigma, \mu)$ is
isometrically isomorphic to another $L_1$-space.

In 1970 Wulbert \cite{W70} extended Theorem~\ref{ando} for arbitrary measure spaces
under the additional assumption that projection $P$ is contractive both in
$L_{p}$-norm and $L_{\infty}$-norm.  At the same time, Tzafriri \cite{T69}
extended Theorem~\ref{isometryLp} to $L_{p}$-spaces on arbitrary
measure spaces.  

Dinculeanu and M.M. Rao studied generalizations of
Theorems~\ref{condexp} and \ref{ando} to $L_{p}(\Omega, \Sigma, \mu)$,
 where
$\Sigma$ is not assumed to be a $\sigma$-algebra but only a $\delta$-ring,
i.e. $\Sigma$ is closed under finite unions, differences and countable
intersections and $\mu$ is a finitely additive general measure (not
necessarily finite).  Dinculeanu \cite{D71} showed that Theorem~\ref
{condexp} still holds in this more general setting, and Dinculeanu and Rao
\cite{DR72} obtained analogues of Theorem~\ref{ando} under additional
assumptions that projection $P$ besides being contractive is also positive,
or the range of $P$ equals to $L^{p}(\Lambda)$ for some sub-$\delta$-ring
$\Lambda \subset \Sigma$, or $P$ satisfies the  averaging identity $(6)$.

In 1974, Bernau and Lacey \cite{BL74} (see also 
\cite{BLcharLp,Lacey})
gave new unified
presentation of the proof of Theorem~\ref{ando} and 
Theorem~\ref{isometryLp} generalized to arbitrary measure spaces.  
When $1<p<\infty$, their approach does not rely on the reduction to the
case of $L_{1}$, but instead they use duality arguments depending on
smoothness of $L_{p}$ and $L_{q}=(L_{p})^{*}$.  Their proof uses the
following two key lemmas:

\begin{lem}
\lb{norming}
\cite[Lemma~1]{Ando},
\cite[Lemma~2.2]{BL74}
Suppose $1<p<\infty$ and let $P$ be a contractive projection on
$L_{p}(\Omega, \Sigma, \mu)$.  Then $f \in \calR(P)$ if and only 
if $J(f) \in
\calR(P^{*})$ (here $J$ denotes the duality map in $L_{p}$, see
Definition~\ref{defJ}).
\end{lem}

\begin{pf*}{Sketch of Proof}
Let $f \in \calR(P)$.
Notice that we have
\begin{align*}
\Vert f \Vert_{p}^{2} &= \left< f,J (f) \right> = \left< Pf,
J(f) \right> \\
&= \left< f, P^{*}(J(f) \right> \leq \Vert f \Vert_{p} \cdot
\Vert P^{*} (J(f)) \Vert_{q} \\
&\leq \Vert f \Vert_{p} \cdot \Vert J (f) \Vert_{q} = \Vert f \Vert^{2}_{p}\ .
\end{align*}

Since $L_{p}$, $1<p<\infty$, is smooth we conclude that
$$
P^{*}(J(f)) = J(f)
$$
i.e.
$$
J(f) \in \calR(P^{*}).
$$

Since $L_{q}=(L_{p})^{*}$ is also smooth we similarly obtain that $J(f) \in
\calR(P^{*})$ implies $J(J(f))=f \in \calR(P^{**}) = \calR(P)$.

Notice that this proof is valid not only in $L_{p}$, but in any smooth,
reflexive Banach space $X$ with a smooth dual.
\end{pf*}

\begin{lem}
\cite[Lemma~2.3(i)]{BL74}
\lb{signum}
Suppose $1<p<\infty, p \neq 2$. Let $P$ be a contractive projection on
$L_{p}(\Omega, \Sigma, \mu)$.  If $f,g \in \calR(P)$ then $\vert f \vert 
\sgn g
\in \calR(P)$.
\end{lem}

\begin{pf*}{Idea of proof}
The proof of this important lemma is somewhat technical.  It involves an
inductive procedure as follows:

Set $k_{0} =g$.
Then, by Lemma~\ref{norming},
$J(f), J(f+ \lambda k_{0}) \in R(P^{*})$ for every $\lambda \in \Bbb R$.
It is then shown that
$$
\lim \limits_{\lambda \to 0} \frac{1}{\lambda} (J(f + \lambda k_{0})- J(f))
$$
exists a.e. and in $L_{q}$.

Set
\begin{align*}
g_{0} &= \lim \limits_{\lambda \to 0} \frac{1}{\lambda} (J(f + \lambda
k_{0})-J(f)) \\
k_{1} &= J^{-1}(g_{0})
\end{align*}

Then $g_{0}\in R(P^{*})$ and, by Lemma~\ref{norming}, $k_{1} \in \calR(P)$.
Define inductively
$$
k_{n+1} =J^{-1}(\lim \limits_{x \to 0} \frac{1}{\lambda}(J(f + \lambda
k_{n})-J(f))).
$$

Then Bernau and Lacey show that $\lim \limits_{n \to \infty} k_{n}$ exists
in $L_{p}$ and that
$$
\lim \limits_{n \to \infty} k_{n} = \vert f \vert \sgn g.
$$

Since $k_{n} \in \calR(P)$ for each $n$, the lemma is proven.
\end{pf*}

Using Lemma~\ref{signum} it is not difficult to deduce:

\begin{lem}
\cite[Lemma~3.1]{BL74} \lb{subalgebra}
Suppose $1<p<\infty$, $p \neq 2$ and let $P$ be a contractive projection on
$L_{p}(\Omega, \Sigma, \mu)$.  Define $\Sigma_{0} = \{\supp f: f \in \calR(P)\}
\subset \Sigma $.  Then $\Sigma_{0}$ is a $\sigma$-subalgebra of $\Sigma$.
\end{lem}

In the next step of the proof of Theorem~\ref{ando} Bernau and Lacey show that
there exists a function $h$ in $\calR(P)$ with maximal support and that 
$\calR(P)$ is
isometrically isomorphic to $L_{p} (\supp h, \Sigma_{0}, \vert h \vert^{p}
\mu_{\Sigma_{0}})$, which leads them to the final conclusion.

\subsection{Results coming from the approximation theory and the nonlinear approach}
Contractive linear projections in $L_p$ were also studied from the point of
view of nonlinear analysis.  This comes from the fact that the existence of
a contractive linear projection onto a subspace is intrinsically related to
the metric projection onto a complementary subspace. 
In Section~\ref{metric} below we discuss this relationship in greater
detail.  Here we just note the fact that was used (and observed) by de
Figueiredo and Karlovitz \cite{FK}:
 
\begin{prop} \lb{proximity1}
Let $X$ be a normed space and $P$ be a linear projection on $X$ with $\codim \calR
(P)=1$.  Then $P$ has norm one if and only if $I-P$ is a metric projection
onto $\Ker(P)$, i.e. for each $x \in X$ and $y \in \Ker (P)$, $\|x- (I-P)
x \| \leq \| x-y \|$ (see Definition~\ref{defprox}).

Thus, if $P$ is a norm one projection, then there exists a linear selection 
of a metric projection onto $Ker(P)$.
\end{prop}

De Figueiredo and Karlovitz used this to obtain the following:

\begin{prop}
\cite[Proposition~3]{FK}
(see also \cite[Proposition~1]{BM77})
\lb{hyperplanesLp}
Let $1< p<\infty,\  p \neq 2$, and $(\Omega,\Sigma, \mu)$ be a $\sigma$-finite
nonatomic measure space.  Then no subspace of $L_p(\Omega, \Sigma, \mu)$ of
codimension one is the range of a linear projection of norm $1$ in
$L_p(\Omega, \Sigma, \mu)$.
\end{prop}

This proposition, of course, also follows quickly from Theorem~\ref{ando}
(also for $p=1$) but the proof in \cite{FK} is significantly simpler.  The
proof in \cite{BM77} also  is an application of
Proposition~\ref{proximity1}. In fact Beauzamy and Maurey obtained a
somewhat stronger result not limited to nonatomic measures:

\begin{th}
\cite[Proposition~1]{BM77} (cf. also \cite[Theorem~5.2]{Deutsch82})
\lb{BM}
Let $1<p<\infty,\  p \neq2$, and let 
$f \in L_q (\Omega, \Sigma, \mu)\ f \neq0
\ (\frac{1}{p} + \frac{1}{q} =1)$.  Then the hyperplane $f^{-1}(0)$ is the
range of a norm one linear projection in $L_p(\Omega, \Sigma, \mu)$ if and
only if $f$ is of the form $f= \alpha \chi_A + \be \chi_B$, where $A$ and
$B$ are atoms of $\mu$ and $ \alpha, \be$ are scalars.
\end{th}

Proposition~\ref{hyperplanesLp} for $p=1$ was proven   using
the methods of $L-$summands \cite[Corollary~IV.1.15]{HWW}.

Contractive projections in $L_{p}$ were used to construct monotone and
unconditional bases in $L_{p}$.  Such constructions depend on the study of
increasing sequences of contractive projections (we say that a sequence of
projections $\{P_{j}\}^{\infty}_{j=1}$ is increasing if
$P_{i}P_{j} = P_{\min\{i,j \}}$ for
all $i,j \in \Bbb N$), see \cite{DO75,PR75,Burk84,Doust89} and the survey
\cite{Doust88} for the discussion of different unconditionality properties
for contractive projections.  We are not aware of any generalizations of
these constructions to spaces other than $L_{p}$, maybe because
contractive projections are not well understood in spaces other than
$L_{p}$ (see Section~\ref{function}).

\subsection{Bibliographical remarks}        
There are many papers that deal with different properties of contractive
projections in Lebesgue spaces $L_p$ in the situation when the projection
satisfies some additional conditions.  We do not have the space here to
describe all the relevant results, we just list the papers that the author
of this survey is aware of:
\cite{BCS72,BS72,B74,BL77,BL85,Sp86,St86,Popov87}.

\section{Nonatomic K\"othe Function Spaces}
\lb{function}

To this day there are no classes of K\"othe function spaces other than
$L_{p}$ where the form of contractive projections is fully characterized.
In this section we present the history of different partial results
describing different properties of norm one projections and
$1$-complemented subspaces of different K\"othe function spaces.

\subsection{Results related to the ergodic theory}
The first general results about existence of norm-one projections are the
consequence of the classical Mean Ergodic Theorem (see
\cite[Section~VIII.5]{DS1}, this theorem goes back to 1930's, see the
excellent bibliographical notes in \cite[Section~ VIII.10, p. 728]{DS1}).
We have:

\begin{th}
\lb{ergodic}
Let $X$ be a reflexive Banach space and $T$ be a norm one operator on $X$.
Then the operators
$$
A_{T,n}=\frac{1}{n+1} \sum \limits^{n}_{j=0} T^j
$$
converge strongly to a norm one projection onto the space $F_T = \{x \in X:
T(x)=x\}$.
\end{th}

Also Lorch \cite{L39}  studied monotone sequences of projections
$\{P_{n}\}^{\infty}_{n=1}$ whose norm has a common bound
(i.e. $\exists K \in \Bbb R$ such that $\Vert P_{n} \Vert \leq K$
for all $n \in \Bbb N$).

\begin{defn}
We say that a sequence of projections $\{P_{n}\}^{\infty}_{n=1}$ is {\it
increasing} (resp. {\it decreasing}) if for all $n,m \in \Bbb N$
$P_{n}P_{m} = P_{\min(n,m)}$ (resp. $P_{n}P_{m}=P_{\max(n,m)})$.
\end{defn}

Lorch showed in particular that:

\begin{th}
\lb{lorch}
Suppose that $X$ is a reflexive Banach space.  Let
$\{P_{n}\}^{\infty}_{n=1}$ be a monotone sequence of contractive
projections.  Then

\begin{itemize}
\item[(a)]
If the  sequence $\{P_{n}\}^{\infty}_{n=1}$ is increasing then $Y = \overline
{\bigcup \limits_{n=1}^{\infty} \calR (P_{n})}$ is $1$-complemented in $X$.
\item[(b)] If the sequence $\{P_{n}\}_{n=1}^{\infty}$ is decreasing than $Y =
\bigcap \limits_{n=1}^{\infty} \calR (P_{n})$ is $1$-complemented in $X$.
\end{itemize}
\end{th}

(Recall that $\calR (P))$ denotes the range of the projection $P$).

The study of contractive projections in connection with the ergodic theory
is still very active, it includes in particular a study of Korovkin sets
and Korovkin approximations, see e.g. \cite{B74,W75,CFS79,Re00} and their
references (\cite[Section~VI]{CFS79} contains a short overview of related 
results).

\subsection{Other general results}
Next result  illustrates that $1$-complemented
subspaces can be very rare and unstable (see also Theorem~\ref{Bh2} and
Section~\ref{nonexistence}).  Lindenstrauss \cite{L64} showed the
following striking fact:

\begin{th}
\lb{L64}
There exist Banach spaces $Z \supset X$ with $\dim Z / X = 2$ such
that
\begin{itemize}
\item[(1)]
For every $\e > 0$ there is a projection with norm $\leq 1 +
\e$ from $Z$ onto $X$,
\item[(2)]
For every $Y$ with $Z  \supsetneq Y  \supsetneq X$ there is a projection
of norm one from $Y$ onto $X$,
\item[(3)]
There is no projection of norm $1$ from $Z$ onto $X$.
\end{itemize}
\end{th}

Thus there does not exist any limiting process that can be applied to
obtain $1$-comp\-le\-ment\-ed subspaces out of 
$(1 + \e)$-comp\-le\-ment\-ed
subspaces.

Results similar to Theorem~\ref{lorch}
were considered in the context of nonlinear projections and with emphasis
on what kind of sequences of (both linear and nonlinear) contractive
projections and products of contractive projections converge to a
contractive projection, see \cite{H62} for a study of products of
projections in Hilbert space, and \cite{BR77,Reich83,DKR91} for the initial
results on products of contractive projections in reflexive Banach spaces (see
also the survey \cite{DR92}).
The literature on this subject continues to grow but the discussion of the
results in this direction is beyond the scope of the present survey.

Next result that we want to mention here is the following fact observed by
Cohen and Sullivan in 1970 \cite{CS70}:

\begin{th}
\lb{uniqueness}
A subspace of a smooth space $X$ can be the range of at most one projection
of norm one.
\end{th}

The question of uniqueness have been subsequently studied in spaces which
are not necessarily smooth, see \cite{O74,O75,O77,OL90}.

\subsection{Characterizations using the duality map}
In 1975 Calvert proved a very important characterization of
$1$-complemented subspaces of a reflexive Banach space $X$ with $X$ and
$X^{*}$-strictly convex.  He showed:

\begin{th}
\lb{calvert}
\cite{Cal75}
Suppose that $X$ is a strictly convex reflexive Banach space with strictly
convex dual $X^{*}$.  Let $J: X \to X^{*}$ be the duality map (i.e. $\Vert
Jx \Vert = \Vert x \Vert, \left< Jx,x \right> = \Vert x \Vert^{2}$, 
see Definition~\ref{defJ}).  Then
a closed linear subspace $M$ of $X$ is the range of a linear contractive
projection if and only if $J(M)$ is a linear subspace of $X^{*}$.
\end{th}

Implication $``\Longrightarrow$'' in Theorem~\ref{calvert} is fairly
immediate.  Indeed, this really is the statement of Lemma~\ref{norming}
which says that if $M=\calR (P)$, where $P$ is a contractive projection,
then $J(M)= \calR (P^{*})$ and thus $J(M)$ is a linear subspace of $X^{*}$.
It is easy to check that the proof of Lemma~\ref{norming} is valid not only
in $L_{p}\ (1<p<\infty)$, but in any smooth reflexive Banach space with
smooth dual (cf. also \cite[Theorem~8]{CS70}).

Implication $``\Longleftarrow$'' follows from the following generalization of
Proposition~\ref{proximity1}
(cf. Section~\ref{metric}):

\begin{prop} \lb{proximity2}
\cite{Cal75}, \cite[Lemma~8]{BM77}
Suppose that $X$ is reflexive, smooth and strictly convex.  Let $P$ be a
projection on $X$.  Then $I-P$ has norm $1$ if and only if $P$ is the
metric projection
(cf. Definition~\ref{defprox}).
\end{prop}

Theorem~\ref{calvert} proved extremely useful in describing the form of
$1$-complemented subspaces in various Banach spaces.  We already saw that
the proof of Bernau and Lacey in $L_{p}$ was based on one part of this
theorem (the key Lemma~\ref{norming}) and it will be used frequently in the
results described in Section~\ref{sequence}.  In fact Theorem~\ref{calvert}
 is so
important that it was proved again (with different methods) by Arazy and
Friedman as a starting point for their study of contractive
projections in $\calC_{p}$ \cite{AF92}.

In 1977 Calvert extended Theorem~\ref{calvert} to general Banach spaces
without the assumption of reflexivity or strict convexity.  He proved:

\begin{th}
\cite{Cal77}\lb{calvert2}
Let $X$ be a Banach space over $\Bbb C$ or $\Bbb R$.  Let $M$ be a closed
linear subspace of $X$.  $M$ is the range a contractive projection if and
only if there exists a $weak^*$-closed linear subspace $L$ of $X^*$ with $M
\subset J^{-1} (L)$ and $L \subset \overline{J(M)}^{\| \cdot \|_{X^{*}}}$.
\end{th}

This appears to be the most general characterization of contractively
complemented subspaces known.

\subsection{Characterizations obtained from the nonlinear approach}
Contractive projections were also studied in the nonlinear theory of Banach
spaces.  We do not have the necessary space here to describe all the
interesting results concerning nonlinear contractive projections (see e.g.
the survey \cite{P83} and the excellent exposition in \cite{Sine87}) but 
we do
want to mention a couple of results concerning linear projections which
were proven as a ``side bonus'' of the nonlinear approach. 

We start from the collection of facts when the existence of a nonlinear 
contractive projection from $X$ onto a linear subspace $Y\subset X$ implies
the existence of a linear 
contractive projection from $X$  onto $Y$ (these facts are listed in
\cite[Proposition~2.1]{Sine87} and are described as ``folklore'' or
``should be folklore'').

\begin{prop}
Let $X$ be a Banach space and $Y$ be a closed linear subspace of $X$.
\begin{itemize}
\item[$(a)$] If $\codim Y=1$ and there exists a nonlinear 
contractive projection from $X$ onto  $Y$ then $Y$ is 
1-complemented (linearly) in $X$.
\item[$(b)$] If $X$ is smooth and if $P:X\lra Y$ is a contractive retraction
onto $Y$ then $P$ is a linear contractive projection.
\end{itemize}
\end{prop}

We also have the following deep fact:
\begin{th}\cite{L64a}
Let $X$ be a Banach space and $Y$ be a closed linear subspace of $X$
such that there exists a nonlinear 
contractive projection from $X$ onto  $Y$. If $Y$ is a conjugate space
 then $Y$ is 
1-complemented (linearly) in $X$.
\end{th}

 Beauzamy
\cite{B75m} introduced the following notions:

\begin{defn}
Let $M$ be a subset of a Banach space $X$.  A point $x \in X$ is called
{\it minimal with respect to $M$} if for all $y \in X\setminus \{x\}$ there exists
at least one $m \in M$ with $\Vert m-y \Vert > \Vert m-x \Vert$.

A set of all points minimal with respect to $M$ is denoted $\min(M)$.
Clearly $M \subset \min(M)$.  A set $M$ is called {\it optimal} if $M=
\min(M)$.
\end{defn}

Beauzamy and Maurey proved that this notion is closely related to
contractive projections.

\begin{th} (\cite{BM77}, see also \cite{B75} for the $``\Longrightarrow$''
direction) \lb{optimal}
Suppose that $X$ is reflexive, strictly convex and smooth.  Then the closed
subspace $Y \subset X$ is optimal if and only if $Y$ is the range of a
linear norm one projection on $X$.
\end{th}

The proof uses Proposition~\ref{proximity2}.

Theorem~\ref{optimal}
has been extended by Godini who weakened the conditions on $X$, but
considered a slightly more restrictive definition of minimal sets.

\begin{defn}
Let $X$ be a real normed linear space and $Y$ a linear subspace of $X$.  To
each nonempty subset $M \subset Y$ we assign a subset $M_{Y,X} \subset X$
defined as follows:

$M_{Y,X} = \{x \in X :$ for all $y \in Y\setminus \{x\}$ there exists $m \in M$
such that $\Vert y - m \Vert > \Vert x - m \Vert\}$.

Thus the set $\min(M)$ defined by Beauzamy and Maurey equals $M_{X,X}$.
\end{defn}

Godini proved that

\begin{th} \cite[Theorem~2]{G80} \lb{godini}
Let $X$ be a normed linear space and $Y$ a closed linear subspace of $X$.
A necessary, and if every element of $Y$ is smooth in $X$ also sufficient
condition for the existence of a norm one linear projection $P$ of $X$ onto
$Y$ is that $Y_{Y,X} = Y$.
If every element of $Y$ is smooth in $X$, then there exists at most one norm
one linear projection of $X$ onto $Y$.
\end{th}

Godini used this result to characterize spaces $X$ which are
1-complemented in $X^{**}$, provided every element of $X$ is smooth in
$X^{**}$.

\subsection{Relations with metric projections} \lb{metric}

As we mentioned a few times above there is a fruitful line of investigating
contractive linear projections in connection with metric projections.  Both
norm-one projections and metric projections (see Definition~\ref{defprox})
are natural generalizations of orthogonal projections from Hilbert spaces
to general Banach spaces.  Thus, not surprisingly, there is an intrinsic
connection between them.  We have the following very clear 
but important fact:

\begin{prop}
\lb{proximity}
Let $X$ be a normed linear space and $M$ a linear subspace of $X$.  Then
for all $x \in X$ and any $y \in P_M(x)$
$$
\|x-y\| \leq \|x\|
$$
and
$$
\| y \| \leq 2 \|x \|.
$$
\end{prop}

In this proposition we do not assume that $P_M(x)$ has a linear selection
(see Definition~\ref{defsel}).  The earliest explicit reference that 
we could find
for Proposition~\ref{proximity} is \cite{S67} (cf.
\cite[Theorem~4.1]{S74}).  However different (usually weaker) versions of
Proposition~\ref{proximity} were observed on many occasions by different
authors, see Propositions~\ref{proximity1}, \ref{proximity2} above and the
excellent exposition of this topic with full references in \cite{S74} (cf.
also \cite{Deutsch83}).  For us the most important is the following
corollary of Proposition~\ref{proximity}.

\begin{prop}
\lb{contractive}
Let $X$ be a normed linear space.  Let $P$ be a linear projection on $X$.
Then $\| P \| =1$ if and only if $I-P$ is a selection of a metric
projection onto $\Ker P \subset X$. (here $I$ denotes the identity operator
on $X$).

In particular, if $\|P \| =1$ then there exists a linear selection of a
metric projection onto $\Ker P$.
\end{prop}

We include the simple proof below:
\begin{pf}
(cf. e.g. \cite{FK})
Suppose that $P$ is a linear projection on $X$ with $\| P \|=1$.  Then for
each $x \in X$ and $y \in \Ker P$
$$
\|x - (I-P)x \| = \| Px \| = \| P(x-y) \| \leq \| x-y \|.
$$
Thus $\Ker P$ is a proximinal subspace of $X$ (cf.
Definition~\ref{defprox}) and $I-P$ is a selection of a metric projection
onto $\Ker P$.

Conversely, if $I-P$ is a selection of a metric projection onto $\Ker P$
then for each $x \in X$
$$
\| Px \| = \| x-(I-P)x \| \leq \| x -0 \| = \| x \|.
$$

Hence $\| P \| =1$.

The proof is finished by noting that if $P$ is linear so is $I-P$.
\end{pf}

As we illustrated above Proposition~\ref{contractive} is a very important
tool in obtaining characterizations of contractive projections
(Proposition~\ref{hyperplanesLp}, Theorems~\ref{BM}, \ref{calvert}, 
\ref{optimal}, \ref{godini}, 
\ref{hyperplanes}, \ref{FS}, \ref{finitecodim}).
Proposition~\ref{contractive} was also used to obtain characterizations of
subspaces of $L_p$ which admit  a linear selection of a metric projection
through the analysis of known results about norm one projections in
$L_p$-spaces \cite{Lin85,Song91}.

\subsection{Relations with a notion of orthogonality}
\lb{ortho}

Contractive projections may be treated as a generalization of orthogonal
projections from Hilbert spaces to general Banach spaces.  This has been
explored by Papini \cite{Pap74}, Faulkner, Huneycutt \cite{FH78}, Campbell,
Faulkner, Sine \cite{CFS79}
and  Kinnunen \cite{Kin84}, who considered the following extension
of orthogonality to general Banach spaces:

\begin{defn} \cite[\S 4]{Amir}
Let $X$ be a real Banach space and $x, y \in X$.  We say that $x$ is {\it
orthogonal to} $y$ in the sense of Birkhoff-James (or simply {\it $x$ is
$BJ$-orthogonal to $y$}), denoted by $ x \perp y$ if
$$
\| x \| \leq \|x + \lambda y \|
$$
for all $\lambda \in \Bbb R$.
\end{defn}

The notion of $BJ$-orthogonality was introduced by Birkhoff \cite{B35} and
developed by James \cite{J45,J47a,J47}, cf. also \cite{Amir}.

In general $ x\perp y$ does not imply $y \perp x$, thus we introduce two
notions of orthogonal projections:

\begin{defn}
A projection $P$ on $X$ is called
{\it left-orthogonal}  (resp. {\it right-orthogonal}) if for each $x \in X$
$$
Px \perp (x - Px)
$$
$$\text{(resp.}\ \ (x-Px) \perp Px ).$$
\end{defn}

Papini \cite{Pap74} obtained characterizations of Hilbert spaces among general
Banach spaces using these notions (see Theorem~\ref{Hspace}).

In real Banach spaces we have the following:
\begin{th}
\cite{FH78,Kin84}
\lb{Kin}
Let $X$ be a real Banach space and $M$ a complemented linear subspace of
$X$.
Let $P$ be a linear projection from $X$ onto $M$.
Then:
\begin{itemize}
\item[(a)]
$\| P \| = 1$ if and only if $P$ is left-orthogonal.
\item[(b)]
$P$ is (a selection of) a metric projection if and only if $P$ is
right-orthogonal.
\item[(c)]
If $\| P \| =1$ then $M \perp \Ker P$.
\item[(d)]
If $P$ is (a selection of) a metric projection then $\Ker P \perp M$.
\item[(e)]
If $M \perp \Ker P$ and $X = M \oplus \Ker P$ then $\| P \| = 1$.
\item[(f)]
If $\Ker P \perp M$ and $X = M \oplus \Ker P$ then $P$ is (a selection of) a
metric projection.
\end{itemize}
\end{th}

In particular, as a corollary Kinnunen obtained another proof of
Proposition~\ref{contractive}.  Also, as an application of
Theorem~\ref{Kin}, he obtained a characterization of norm one projections
of a finite rank.
For this we
 need the following definition:

\begin{defn}
\cite{Singer1}
Let $\{x_n\}$ be a basis of a Banach space $X$.  Then the {\it sequence of
coefficient functionals $\{f_n\}$ associated to the basis $\{x_n\}$} is
defined by
$$
f_j(x_k) = \delta_{jk}
$$
for all $j, k$, where $\delta_{jk}$ denotes the Kronecker delta.

A basis $\{x_n \}$ is called {\it normal} if $\| x_n \|_X = \| f_n
\|_{X^{*}}=1$ for all $n$.
\end{defn}

By \cite[Theorem~II.2.1]{Singer1} a basis $\{x_n\}$ is normal if and only
if $\| x_n \| = 1$ and $x_n \perp \overline{span} \{ x_1, \cdots,
x_{n-1}, x_{n+1}, \cdots \}$ for each $n$.

Further, by \cite[Theorem~II.2.2]{Singer1} every finite dimensional Banach
space has a normal basis.  Kinnunen proved the following characterization
of norm one projections in terms of normal bases:

\begin{th}
\cite[Theorem~5.3]{Kin84}
Let $X$ be a real Banach space, $M$ be a $1$-complemented linear subspace
of $X$ with $\dim M = n < \infty$, and $P: X \ONTO M$ be a norm one
projection.  Then $P$ is of the form:
$$
Px = \sum \limits^n_{k=1} f_k (x) u_k
$$
for all $x \in X$, where $\{u_k\}^n_{k=1}$ is a normal basis of $M$ and
$\{f_k \}^n_{k=1} \subset X^*$ are norm one functionals such that
$f_k(u_i)= \delta_{ki}$ for all $k, i = 1, \cdots, n$.
\end{th}

The problem of describing the basis structure of $1$-complemented subspaces
in general real Banach spaces is still open, see Section~\ref{sequence} for
the discussion of known results in both real and complex Banach spaces.

\subsection{Relations with isometries} \lb{seciso}

Here we consider the following question:

\begin{question}\lb{qisometry}
Let $X$ be a Banach space and $T\INTO X$ be an isometry. Is the range
of $T$, $Y=T(X)$, 1-complemented in $X$?
\end{question}

This question has an affirmative answer in Hilbert spaces 
(see Proposition~\ref{isometryH}) and in $L_p, \ 1\le p<\infty$ 
(see Theorem~\ref{andolattice}(B3)).

Question~\ref{qisometry} for reflexive Banach spaces was posed by Faulkner and
Huneycutt \cite{FH78}. It is known that in $C[0,1]$ a range of an isometry
does not have to be even complemented \cite{Ditor73}.

Question~\ref{qisometry} was considered by Campbell, Faulkner and Sine 
\cite{CFS79} who proved:

\begin{th}
Let $X$ be a reflexive Banach space and $T$ be an isometry on $X$. If
 the range
of $T$, $Y=T(X)$, is 1-complemented in $X$ then $T$ is a Wold isometry.
\end{th}

Here ``Wold isometry'' is defined as follows:

\begin{defn}
Let $T$ be an injective linear map on a Banach space $X$. Then $T$ is called
a {\it unilateral shift} provided there exists a subspace $L$ of $X$ for which
$$X=\bigoplus_{n=0}^\infty T^n(L)\ .$$

An isometry $T$ on $X$ is called a {\it Wold isometry} provided 
$X=M_\infty \oplus N_\infty$, where $M_\infty=\bigcap_{n=1}^\infty T^n(X)$
and $N_\infty=\sum_{n=0}^\infty \oplus T^n(L)$, where $L$ is a complement
for the range of $T$, $T(X)$, in $X$.
 
Then  $T\big|N_\infty$ is a shift and  $T\big|M_\infty$ is a surjective
isometry (sometimes referred to as a unitary operator).
\end{defn}

This definition was introduced in \cite{FH78} and used to study extensions
to reflexive Banach spaces
of the Wold Decomposition Theorem, which says that every isometry on a Hilbert
space is the sum of a unitary operator and copies of the unilateral shift.

Campbell, Faulkner and Sine 
\cite{CFS79}  also gave an example of a $C(K)$ space and a Wold isometry $T$
on $C(K)$ such that the range of $T$ is 2-complemented but not
1-complemented in $C(K)$ (in this example the range of $T$ is even 
finite-codimensional in $C(K)$).

It is not known whether every isometry in a reflexive Banach space
is a Wold isometry and thus Question~\ref{qisometry} is still open (see
also Remark in Section~\ref{subseccharlp}).

\subsection{Characterizations in terms of the conditional 
expectation operators} \lb{sectioncondexp}

Next came a series of results which studied the form of a contractive
projection in terms of the conditional expectation operators.  We have the
following generalization of Theorem~\ref{condexp} from $L_p$ to general
rearrangement invariant spaces.

\begin{th}
\cite[Theorem~2.a.4]{LT2}
\lb{LT}
Let $X$ be a rearrangement invariant function space (see
Definition~\ref{defri}) on the interval $I$, where $I=[0,1]$ or $I=[0,
\infty)$.  Then, for every $\sigma$-subalgebra $\Sigma_0$ of measurable
subsets of $I$ so that the Lebesgue measure restricted to $\Sigma_0$ is
$\sigma$-finite, the conditional expectation operator
$\calE^{\Sigma_{0}}$ is a projection of norm one from $X$ onto the
subspace $X_{\Sigma_{0}}$ of $X$ consisting of all the
$\Sigma_0$-measurable functions in it.
\end{th}

As Lindenstrauss and Tzafriri point out, Theorem~\ref{LT} is a consequence
of Theorem~\ref{condexp} and general interpolation theorem (although they
do give a direct proof of it), and thus it is really valid in any
interpolation space between $L_1$ and $L_\infty$.

\begin{rem}
The statement of Theorem~\ref{LT} appears explicitly in
\cite[Section~11.2]{Du71}, we do not know whether or not this is the first
reference for it.
\end{rem}

Similarly as in the case of $L_p$, a lot of effort has been put into
proving the converse of Theorem~\ref{LT}, i.e. into proving the following
conjecture, which would generalize Theorem~\ref{ando} for $L_p$:

\begin{con}
\lb{condexpcon}
Let $X$ be a rearrangement invariant function space on the interval
$[0,1]$ (so $X$ contains constant functions).  Suppose that $P$ is a
contractive projection on $X$ with $P(\bold 1) =\bold 1$.  Then there
exists a $\sigma$-algebra $\Sigma_0$ of measurable subsets of $[0,1]$ so 
that
$P$ is the conditional expectation operator $\calE^{\Sigma_{0}}$.
\end{con}

In fact many extend Conjecture~\ref{condexpcon} to any function space $X$
for which Theorem~\ref{LT} is valid, i.e. to spaces $X$ with ${\bold 1} \in
X$ and where conditional expectation operators are norm one projections
(Bru and Heinich \cite{BH85} call spaces $X$ with this property {\it
invariant under conditioning}).

Duplissey \cite[Theorem~II.2.1]{Du71} showed that arbitrary conditional
expectation operators are contractive in a K\"othe space $X$ if and only if
all conditional expectation operators with finite dimensional range are
contractive on $X$.  Duplissey also studied Conjecture~\ref{condexpcon} but
with the additional assumption that $P$ is contractive in $L_{\infty}$-norm as
well as in norm of $X$.  He proved:

\begin{th}
\cite[Theorem~II.5.5]{Du71}
\lb{Du}
Let $X$ be a strictly monotone K\"othe function space on a $\sigma$-finite
measure space $(\Omega, \Sigma, \mu)$.
Then the following are equivalent:
\begin{itemize}
\item[(a)]
If $\mu(\Omega)$ is not finite then for each $A \in \Sigma$ such that $0
\neq \chi_A \in X$ and $\mu(A)<\infty$, the element $\chi_A$ is   smooth 
in $X$ and
$$
J(\chi_A) = \frac{{\| \chi_{A} \|}_X}{\mu(A)} \cdot \chi_A .
$$
If $\mu(\Omega)$ is finite, then $\bold {1}$ is a smooth point of $X$ and
$J(\bold 1) =({\| {\bold 1 }\|}_{X}/{\mu(\Omega)}) \cdot \bold 1$.
\item[(b)]
Each positive contractive projection on $X$ such that $\|Pf \|_\infty \leq
\| f \|_\infty$ for all $f \in X$ is a conditional expectation operator.
\end{itemize}
\end{th}

First developments in the study of Conjecture~\ref{condexpcon} without any
additional assumptions on a contractive projections are due to Bru and
Heinich, \cite{BH85,BH89} and Bru and Heinich and Lootgieter \cite{BHL86}.
We will outline here the results in \cite{BH89} which contain and expand
on earlier work \cite{BH85,BHL86}.  The authors start from generalizing the
crucial Lemma~\ref{signum}, which was used in the proof of Bernau and Lacey
in $L_p$.  They prove:

\begin{prop}
\cite[Proposition~8]{BH89} (cf. \cite[Proposition~7]{BH85})
\lb{BH1}
Let $X$ be an order continuous K\"othe function space such that conditional
expectation operators are contractive on $X$.  Assume that norm of $X$ is
twice differentiable at $\bold 1$ and $\| {\bold 1} \|_{X}=  1$.  If $X
\subset X^{*}$ then there exists a constant $k \geq 0$ such that for all $f
\in X$
$$
\lim \limits_{\e \to 0} \frac{1}{\e}(J({\bold 1} + \e f)-
{\bold 1}) =k(f-(\int \limits_I f(t)dt)\cdot {\bold 1})
$$
where the limit is taken in norm of $X^*$, and recall that $J: X \to X^*$
denotes the duality map (see Definition~\ref{defJ}).
\end{prop}

The statement of Proposition~\ref{BH1} captures the most essential element
of the proof of Lemma~\ref{signum}, but it's proof is  much less technical;
it makes an elegant use of differentiability of $\| \cdot \|_X$ and
$J(\cdot)$ at ${\bold 1}$.

To generalize the exact statement of Lemma~\ref{signum}, Bru and Heinich
introduce the following definition:

\begin{defn}
Let $X$ be an order continuous, smooth K\"othe function space with ${\bold
1} \in X$ and $\|{\bold 1} \|_X = 1$.  Then $X$ is called {\it D-concave}
(Bru and Heinich use $D$ to denote the duality map which in this survey
is denoted by $J$, following the most of English-language literature; so it
would be natural for us to use ``$J$-concave'') if
\begin{itemize}
\item[(i)]
$X^* \subset X$, the inclusion map is continuous and
$$
\lim \limits_{c \to \infty} 
\sup\{\|f \cdot \chi_{\{| f |>c\}} \| : f \in X^*,
\| f \|_{X^{*}} =1\}=0
$$
(i.e. the unit ball of $X^*$ is $X$-equiintegrable),
\item[(ii)]
$J(J(f))= f \Longleftrightarrow f = {\sgn(f)}/{\| \sgn (f) \|_{X^*}}$.
\end{itemize}

$X$ is called {\it $D$-convex} if $X^*$ is $D$-concave.
\end{defn}

Notice that $L_p[0,1],\ 1<p<2$ are $D$-concave and $L_p[0,1],\  2<p<\infty$ 
are
$D$-convex.

Now we are ready to state the generalization of Lemma~\ref{signum}.

\begin{prop}
(cf. \cite[Theorem~11]{BH89})
\lb{BHsignum}
Let $X$ be a $D$-concave K\"othe function space, and let $P$ be a
contractive projection on $X$ with $P({\bold 1})={\bold 1}$.  Then if $f
\in \calR(P)$ then $\sgn(f) \in \calR(P)$.
\end{prop}

Then using the analogue of Lemma~\ref{subalgebra} Bru and Heinich obtain the
generalization of Theorem~\ref{ando} for constant preserving contractive
projections:

\begin{th}
\cite[Theorem~13 and its Corollary]{BH89}
\lb{BH2}
Suppose that $X$ is a $D$-concave K\"othe function space such that norm of
$X^*$ is twice differentiable at $\bold 1$, or $X$ is a $D$-convex K\"othe
function space such that norm of $X$ is twice differentiable at $\bold 1$.
Let $P$ be a contractive projection on $X$ with $P({\bold 1})={\bold 1}$.
Then $P$ is a conditional expectation operator.
\end{th}

As a corollary they obtain a characterization of constant preserving
contractive projections is special Orlicz spaces:

\begin{th}
\cite[Proposition~28]{BH89}
\lb{BH3}
Let $\varphi$ be an Orlicz function which is twice differentiable on $\Bbb
R_{+}$ and such that $\varphi'' $ is either strictly increasing to infinity
or $\varphi''$ is strictly decreasing to $0$.  Then every contractive
projection $P$ with $P({\bold 1})={\bold 1}$ on the Orlicz function space with
Luxemburg norm $L_\varphi$, or with Orlicz norm $L_{\varphi,O}$ is a
conditional expectation operator.
\end{th}

They also obtain the same conclusion under the assumption that Orlicz
function $\varphi$ has a continuous strictly increasing derivative
$\varphi^{\prime}$ so that $\varphi^{\prime}$ is of  the concave type
(i.e. there exist constants
$\gamma,\  t_0 > 0$
so that for all
$\lambda,\  0< \lambda \leq 1$
and all
$t \geq t_0  \ \ , \varphi^\prime(\lambda x)/(\lambda x)
\geq \gamma(\varphi^\prime(x)/x)$,
for all
$t,\  0<t<1,\ \varphi^\prime (t)>t$;
for all $t>1,\  \varphi^\prime (t)<t$; \ 
$\varphi^{\prime}$ is differentiable at $t=1$ and
$\lim \limits_{t \to \infty} \varphi^\prime (t)= \infty$)
or the inverse $(\varphi')^{-1}$ is of the concave type
\cite[Theorem~4]{BH85},  
\cite[Application]{BHL86}.

One might say that the restrictions on K\"othe function space $X$ in
Theorem~\ref{BH2} and on the Orlicz function in Theorem~\ref{BH3} are
somewhat severe, however these results are the most general results known
about the form of contractive projection $P$ with $P({\bold 1})={\bold 1}$.
Nothing, outside of $L_p$, is known about contractive projections which do
not satisfy $P({\bold 1})= {\bold 1}$.  In particular it is not known which
functions in $X$ can be the image of ${\bold 1}$.

Next development in the study of contractive projections in terms of
conditional expectation operators is due to P. Dodds, Huijmans and de
Pagter \cite{DHP90} who obtained very general results under the assumption
that contractive projection $P$ is positive or that the range of $P$ is a
sublattice.  They extended Theorem~\ref{P1=1} to very general lattices, but
under the restriction that $P$ is positive.

\begin{th}
\cite[Proposition~4.7]{DHP90}
\lb{P1=1positive}
Let $X$ be a K\"othe function space on a finite measure space $(\Omega,
\Sigma, \mu)$ and let $P:X\to X$ be a linear map.  Then the following
statements are equivalent:
\begin{itemize}
\item[(a)]
There exists a $\sigma$-algebra $\Sigma_{0} \subset \Sigma$ such that $P$
is the conditional expectation operator $\calE^{\Sigma_{0}}$.
\item[(b)]
$P$ is a positive order continuous projection with $P({\bold 1})={\bold 1}$
and $P'({\bold 1})={\bold 1}$ (here $P': X' \to X'$ denotes the dual
mapping to $P$ defined on the K\"othe dual  $X'$).
\end{itemize}
\end{th}

\begin{cor}
\cite[Corollary~4.9]{DHP90}
Let $X$ be a K\"othe function space such that norm on $X$ is smooth at
$\bold 1$ and $\|{\bold 1} \|_{X} \cdot \|{\bold 1} \|_{X'} =
\mu(\Omega)$.  If $P$ is a positive order continuous contractive projection
with $P({\bold 1})= {\bold 1}$, then there exists a $\sigma$-subalgebra
$\Sigma_0 \subset \Sigma$ such that $P = \calE^{\Sigma_{0}}$.
\end{cor}

This corollary is a significant extension of Theorem~\ref{Du}.

Next Dodds, Huijmans and de Pagter obtain characterizations of contractive
projections onto a sublattice.

\begin{th}
\cite [Corollary~4.14]{DHP90}
Let $X$ be a K\"othe function space with order continuous norm and let $P$
be a contractive projection in $X$ such that $\calR(P)$ is a sublattice.
If $\calR(P)$ contains some strictly positive functions and if the norm $X$
is smooth at all such strictly positive functions then $P$ is a weighted
conditional expectation operator i.e. there exists a $\sigma$-subalgebra
$\Sigma_0 \subset \Sigma, \ 0 \leq w \in X'$ and $0<k \in L_1(\Omega, \Sigma,
\mu)$ with
$$
\calE^{\Sigma_{0}}(wk) =\calE^{\Sigma_{0}}(k)= \bold 1
$$
such that
$$
Pf = k \calE^{\Sigma_{0}}(wf)
$$
for all $f \in X$.
\end{th}

\begin{th}
Let $X$ be a K\"othe function space with order continuous norm such that
the norm is smooth at ${\bold 1}$ and $\|{\bold 1}\|_X \| {\bold 1}
\|_{X'} = \mu(\Omega)$.
If $P$ is a contractive projection in $X$ such that $\calR(P)$ is a
sublattice and $ {\bold 1} \in \calR(P)$, then $P$ is a conditional
expectation operator, i.e. there exists a $\sigma$-subalgebra $\Sigma_0
\subset \Sigma$ so that $P=\calE^{\Sigma_{0}}$.
\end{th}

The paper \cite{DHP90} contains also many very interesting results about what
conditions on positive operators $T$ (which are not necessarily projections
nor contractive) assure that they will be conditional expectation
operators, but these results will not be summarized here.  We finish the
account of work in \cite{DHP90} with results which relate when a 
 contractive projection $P$ on a Banach lattice $X$ is
 positive and when the range of $P$ is a  sublattice of $X$.

\begin{th}
Let $X$ be a Banach lattice and $P$ be a contractive projection on $X$
\begin{itemize}
\item[(a)] \cite[Lemma~4.5]{DHP90}, \cite[Theorem~II.3.2(i)]{Du71} 
If $X$ is strictly monotone and $P$ is positive then $\calR(P)$ is a
sublattice.
\item[(b)]
\cite[Remark after Lemma~4.5]{DHP90}
There exists $X$ not strictly monotone (e.g. $X = \ell_\infty^3$) and $P$
positive with $\calR(P)$ not a sublattice.
\item[(c)]
\cite[Proposition~4.10]{DHP90}
If $X$ is smooth and $\calR(P)$ is a sublattice then $P$ is positive.
\item[(d)]
\cite[Example~4.11]{DHP90}
There exists $X$ nonsmooth ($\dim X=3$, ball of $X$ is a dodecahedron, thus
$X$ is also non-strictly monotone, but it is symmetric) and non-positive
$P$ with $\calR(P)$ a sublattice.
\item[(e)]
\cite[Proposition~4.13]{DHP90}
If $X$ is order continuous and $\calR(P)$ is a sublattice such that
$\calR(P)$ contains some strictly positive functions from $X$ and norm of
$X$ is smooth at all strictly positive functions in $\calR(P)$, then $P$ is
positive.
\item[(f)]
\cite[Proposition~4.15]{DHP90}
If $X$ is order continuous and $\calR(P)$ is a sublattice such that there
exists a strictly positive function $w \in \calR(P)$ so that norm of $X$ is
smooth at $w$ and $J(w)$ is strictly positive then $P$ is positive  (this
is satisfied for example if ${\bold 1} \in \calR (P)$ and norm of $X$ is
smooth at $\bold 1$).
\end{itemize}
\end{th}

\subsection{Nonexistence of 1-complemented subspaces of finite codimension}

We finish this section with results about non-existence of contractive
projections onto subspaces of finite codimension in a K\"othe function
space $X$, which extend
Proposition~\ref{hyperplanesLp}.
We have:
\begin{th}
\cite[Theorem~4.3]{KR}
(cf. also \cite[Theorem~2]{pams})
\lb{hyperplanes}
Suppose that $X$ is a separable, real order-continuous K\"othe function
space on $(\Omega, \Sigma, \mu)$, where $\mu$ is nonatomic and finite. Then
the hyperplane $M=f^{-1}(0),(f \in X^*)$ is $1$-complemented if and
only if there exists a nonnegative measurable function $w$ with $\supp w =
B = \supp f$, so that for any $x \in X$ with $\supp x \subset B$
$$
\| x \|_X = ( \int |x|^2 w d \mu)^{\frac{1}{2}},
$$
i.e. there are no $1$-complemented hyperplanes in $X$ unless $X$
contains a band isometric to $L_2$.
\end{th}

Around the same time a similar result was obtained by Franchetti and
Semenov for rearrangement-invariant function spaces, but without the
restriction of separability:

\begin{th}
\cite[Theorem~1]{FS95}
\lb{FS}
Let $X$ be a real rearrangement-invariant function space on $(\Omega,
\Sigma, \mu)$, where $\mu$ is nonatomic and $\mu(\Omega)=1$. Denote by $S$
a rank one projection
$$
Sx = (\int \limits_\Omega x(s) d\mu(s)) {\bold 1}.
$$
Then $\|I-S\|=1$ if and only if $X$ is isometric to $L_2 (\Omega, \Sigma,
\mu)$; i.e. the hyperplane $f^{-1} (0)$, where $f(x)= \int_{\Omega}
x(s)d \mu(s)$ for all $x \in X$, is $1$-complemented in $X$ if and only if
$X=L_2(\Omega, \Sigma, \mu)$.
\end{th}

Next the author of this survey extended Theorem~\ref{hyperplanes} to
subspaces of any finite codimension:

\begin{th}
\cite[Theorem~4]{pams}
\lb{finitecodim}
Suppose $\mu$ is nonatomic and $X$ is a real separable
re\-ar\-ran\-ge\-ment-invariant space on $[0,1]$ not isometric to $L_2$.  Then
there are no $1$-complemented subspaces of any finite codimension in $X$.
\end{th}

The proofs of Theorems~\ref{hyperplanes}, \ref{FS}, \ref{finitecodim} 
all use the
classical Liapunoff Theorem (see e.g. \cite{Rudin}) and facts related to
Proposition~\ref{contractive}.  Moreover, in \cite{KR,pams} the following fact is
used:

\begin{prop} (cf. \cite{KR,R84})\lb{flinn}
Let $X$ be a real Banach space and $P$ be a projection in $X$.  Then $\| I-P
\|=1$  (where $I$ denotes the identity operator) if and only if for all $x
\in X$ there exists $x^* \in X^*$ with $x^* \in J(x)$ and $\left< x^*, Px \right> \geq 0$.  ($X$
does not have to be smooth, see Definition~\ref{defJ}).
\end{prop}

The proof of this fact uses the theory of numerical ranges \cite{BD1,BD2} 
and it
 relates contractive projections to accretive operators in real
Banach spaces. Proposition~\ref{flinn} is
 very useful in characterizing contractive
projections in Banach spaces with $1$-unconditional bases, see
 Section~\ref{sequence}.

\begin{rem}
A careful reader may have noticed that, despite our effort to have as
complete a bibliography as possible, there are several papers concerning
contractive projections in function spaces by M.M. Rao that have not been
referenced here.  The reason for this omission is that, regrettably, many of
results in those papers are not valid in the full generality as stated
there, but the methods are in fact limited only to the $L_{p}-$case, see
\cite{BL74,H80,Lreview}.
\end{rem}

\part{Sequence spaces}

\section{Lebesgue sequence spaces $\ell_p$} \lb{lpsection}

In this section we discuss the development of the study of 
 contractive
projections in the  case of the Lebesgue sequence spaces $\ell_{p}$.

\subsection{General  results} 
First result about 1-complemented subspaces of $\ell_p$ is due to
Bohnenblust  who considered finite-dimensional spaces $\ell_p^n$:

\begin{th} \cite[Theorem~3.2]{Bh41} \lb{Bh}
A subspace $S$ of a n-dimensional space $\ell_p^n$ is 1-com\-ple\-ment\-ed in
$\ell_p^n$ if and only if $S$ is spanned by disjointly supported vectors.
\end{th}

The method of proof of Theorem~\ref{Bh} is technically very complicated, it
involved conditions of Pl\"ucker Grassmann coordinates of the subspace $S$ (we
will not present the definition here).  However the proof, in addition to
Theorem~\ref{Bh}, gives also a characterization of 1-complemented subspaces
of $n$-dimensional subspaces $S \subset \ell_p^n$.  As a corollary
Bohnenblust showed that there exist subspaces of $\ell_p^n$ which do not
have any 1-complemented subspaces (see also Section~\ref{sequence}).

\begin{th} \cite[Theorem 3.3]{Bh41}\lb{Bh2}
Let $1 < p < \infty ,\ p$ not an integer and let $l \in \Bbb N$
be such that $2(2l - 3) < n$.  Then there exists $l$-dimensional
subspaces $S_l$ of $\ell_p^n$ such that only $S_l$ and subspaces of
dimension one are 1-complemented in $S_l$.
\end{th}

The case of infinite dimensional $\ell_p$ is simpler than the case of general
$L_p-$spaces, however the original proofs of Douglas
and Ando do not cover it (as they work only on finite measure spaces).
Subsequent generalizations by Tzafriri and Bernau, Lacey do not consider
this case separately.  The simple proof specifically for $\ell_p$ is
included in \cite[Theorem~2.a.4]{LT1}. We quote the statement of this 
theorem
 below because it
illustrates better than the general case, the geometric properties of
$1$-complemented subspaces which we will try to transfer to other 
  sequence
spaces.

\begin{th}  \cite{LT1} \lb{lp}
Let $1 \leq p< \infty$, $p\neq 2$, and $F \subset \ell_{p}$ be a closed linear subspace of
$\ell_{p}$.  Then the following conditions are equivalent:

\begin{itemize}
\item[(1)] $F$ is $1$-complemented in $\ell_{p}$,

\item[(2)] $F$ is isometric to $\ell_{p}^{\dim F}$,

\item[(3)] There exist vectors 
$\{u_{j}\}_{j=1}^{\dim F}$ of norm one and
the form 
$$u_{j}=\sum_{k \in S_{j}} \lambda_{k} e_{k},$$ 
with 
$S_{j} \subseteq
\Bbb N$, $S_{j} \cap S_{i} = \emptyset$ for ${j \neq i}$,  and 
such that $F =\overline {\span} \{u_{j}\}^{\dim F}_{j=1}$ 
(here $\{e_{k}\}_{k
\in \Bbb N}$ denotes the usual basis of $\ell_{p}$).
\end{itemize}

Moreover, if these conditions are satisfied then the norm one projection
$P: \ell_{p} \ONTO  F$ is given by
$$Px = \sum_{j=1}^{\dim F} u_{j}^{*}(x)u_{j},$$ 
where
$\{u_{j}^{*}\}^{\dim F}_{j=1} \subset X^{*}$ satisfy $\Vert u_{j}^{*} \Vert
= u_{j}^{*} (u_{j})=1$ (i.e. $u_j^*=J(u_j)$).
\end{th}

\begin{pf*}{Sketch of proof}
The proof of Theorem~\ref{lp} $(1) \Longrightarrow  (3)$ presented in
\cite{LT1} follows the essential steps of the proof of Bernau and Lacey for
general $L_{p}$, i.e. it also rests on the key Lemma~\ref{signum}.  However
the finishing step to reach the final conclusion is now much simpler than
in $L_{p}$.

Indeed, let $\Sigma_{0}= \{ \supp f: f \in  F= \calR (P)\} \subset \calP
(\Bbb N)$.  Then Lemma~\ref{signum} implies that if $A,B \in \Sigma_{0}$
then $A \cap B \in \Sigma_{0}$.  Thus, for each $i_{0} \in \Bbb N$ such
that $i_{0}$ belongs to the support of $Pf$, for some $f \in \ell_{p}$,
there is a set $A_{i_{0}} \in \Sigma_{0}$ which is minimal in $\Sigma_{0}$
and such that $i_{0} \in A_{i_{0}}$.  Let $A_{i_{0}} = \supp f_{i_{0}}$
where $f_{i_{0}} \in \calR (P)$ and consider a subspace of $F= \calR (P)$
consisting of all functions $g$ so that $\supp g \subset A_{i_{0}} = \supp
f_{i_{0}}$.  We claim that this subspace of $\calR (P)$ is one dimensional.
Indeed, if it was not one dimensional then there would exist $g \in F =
\calR (P)$ linearly independent with $f_{i_{0}}$ and such that $\supp g
\subset A_{i_{0}}$. For any $i \in A_{i_{0}} \backslash \{i_{0} \}$ we now
can find a linear combination of $g$ and $f_{i_{0}}$ so that $(ag +
bf_{i_{0}})(i)=0$.  Since $ag + bf_{i_{0}} \in F = \calR (P)$ we get $\supp
(ag + bf_{i_{0}}) \in F= \calR (P)$ and $\supp (ag + bf_{i_{0}})
\subseteq A_{i_{0}}$, which contradicts minimality of $A_{i_{0}}$.  Now
let $u_{i_{0}} \in F$ be the unique vector with $\supp u_{i_{0}} =
A_{i_{0}}$ and $\Vert u_{i_{0}} \Vert_{p} =1$.  It is easy to see that (3)
holds.

The final statement of Theorem~\ref{lp} about the form of the projection
$P$ follows from the uniqueness of this projection.
\end{pf*}

The precursor of Theorem~\ref{lp} was proved in 1960
by Pe\l czy\'nski \cite{Pel60}, who showed that in $\ell_{p}$ $ (1<p<\infty)$ the
subspaces that are isometric to $\ell_{p}$ are $1$-complemented in
$\ell_{p}$.

Theorem~\ref{lp} explicitly relates the one-complementability of the
subspace $F$ with the property that $F$ is spanned by disjointly supported
vectors.  In the next section we will analyze such a relation in other
sequence spaces.

\subsection{1-complemented subspaces of finite codimension}

Contractive projections in $\ell_p$ were also investigated from the
approximation theory point of view, as part of the study of minimal
projections.  

\begin{defn}
A projection $P: X \ONTO Y$ is called {\it minimal} if $\| P
\|=\inf \{\| Q \|: Q: X \ONTO  Y, Q \ {\text {projection}}\}$.
\end{defn}

Results on minimal projections appeared in the 1930s, in connection with
  geometry of Banach spaces and they have many applications in numerical
analysis and approximation theory, see the survey \cite{ChP70} for the early
results, and \cite{OL90} for the book length presentation of more modern
developments.
Here we will just give a brief account of results on contractive
projections (which clearly are always minimal).  We start with a result of
Blatter and Cheney, who studied minimal projections onto hyperplanes (i.e.
subspaces of codimension $1$) in $\ell_1$ and $c_0$ \cite{BCh}.  In
particular they proved:

\begin{th}
\cite[Theorem~3]{BCh}
Let $0 \neq f \in \ell_ \infty$.  The hyperplane $Y = f^{-1} (0) \subset
\ell_1$ is a range of a norm one projection in $\ell_1$ if and only if at
most two coordinates of $f$ are different from $0$.  The norm one
projection onto $Y$ is unique if and only if exactly two coordinates of $f$
are different from $0$.
\end{th}

Precisely the same statement is valid for all $p, \ 1\le p <\infty, \ p \ne 2$,
as shown by Beauzamy and Maurey \cite{BM77}, see Theorem~\ref{BM}.

These results have been generalized by Baronti and Papini to subspaces of
arbitrary finite codimension and to arbitrary $p, 1 \leq p < \infty$.  They
proved:

\begin{th}
\cite[Theorem~3.4]{BP88}, 
\cite[Theorem~5.5]{BP91}
\lb{BP}
Let $Y$ be a subspace of $\ell_p\  (1 \leq p < \infty,\  p \neq 2)$ of finite
codimension $\codim Y = n \in \bbN$.  Then $Y$ is $1$-complemented in
$\ell_p$ if and only if $Y$ is the intersection of $n$ $1$-complemented
hyperplanes i.e. if and only if there exist functionals $ f_1, \cdots, f_n
\in (\ell_p)^*$ such that for each $j \leq n$ at most two coordinates of
$f_j$ are different from $0$ and $Y = \bigcap \limits_{j=1}^n f_j^{-1}(0)$.
\end{th}

The proof of this theorem depends on Theorem~\ref{calvert} and is slightly
simpler than the proof of Theorem~\ref{lp} since it is restricted to
subspaces of finite codimension.  Notice that the descriptions given in
Theorems~\ref{lp} and \ref{BP} are equivalent.  This is an intuitively
 straightforward
fact but since we have not seen it in the  literature we present the full
proof below.  Unfortunately, the proof is somewhat techinical. We have:

\begin{prop}
Let $X$ be a Banach  space with basis $\{e_n \}_{n \in \bbN}$ and $Y$ be a
subspace of $X$ with $\codim = n$. Then the following conditions are
equivalent:

\begin{itemize}
\item[(i)]
there exit vectors $\{u_j\}_{j \in \bbN}$ of the form $u_j = \sum
\limits_{k \in S_{j}}  u_{jk} e_k$ with $S_j \subseteq \bbN, S_j \cap S_i
= \emptyset$ for $j \neq i$ and such that $Y= \overline {\span}{ \{u_j\}_{j
\in \bbN}}$.
\item[(ii)]
there exist functionals $f_1, \cdots, f_n \in X^*$ such that
$$
f_j = \sum \limits_{k \in F_{j}} f_{j {k}} e_k^*
$$
with $\card (F_j) \leq 2$ for each $j=1, \cdots, n$ and such that $Y =
\bigcap \limits^n_{j=1} f_j^{-1}(0)$.
\end{itemize}
\end{prop}

\begin{pf}

$(i) \Longrightarrow (ii):$
It is easy to see that since $\{S_j\}_{j \in \bbN}$ are mutually disjoint
$$
\codim Y = \sum \limits^{\infty}_{j=1}(\card (S_j)-1).
$$
Since $\codim Y = n$, all vectors $\{u_j \}_{j \in \bbN}$, except at most
$n$ of them, have a singleton support.

After reordering of $\{u_j \}_{j \in \bbN}$  if necessary, let $m \leq n$ be
such that $ \card (S_j) \geq 2$ for $j \leq m$ and $\card(S_j)=1$ for
$j>m$.

For each $j \leq m$ select $\varphi (j) \in S_j$ and for each $k \in S_j
\backslash \{\varphi (j)\}$ set:
$$
f_k = u_{jk} e^*_{\varphi (j)} - u_{j \varphi (j)} e^*_k
$$

Then:
\begin{equation}
\lb{kernel}
\bigcap \limits_{k \in S_{j} \backslash \{(j)\}} f^{-1}_k (0) = \overline {\span}
(\{u_j\} \cup \{e_i\}_{i \notin S_{j}})
\end{equation}

Indeed, for each $k \in S_j \backslash \{\varphi (j)\}$
$$
f^{-1}_k (0) =\overline {\span}(\{e_i\}_{i \neq k, \varphi(j)} \cup \{ u_{j
\varphi (j)}e_{\varphi (j)} + u_{jk} e_k\})
$$
so $u_j \in f^{-1}_k (0)$ and $\{e_i\}_{i \notin S_{j}} \subset f^{-1}_k
(0)$ for each $k \in S_j \backslash \{\varphi (j)\}$.
By the equality of codimensions we obtain
\eqref{kernel}.

Further
\begin{align*}
Y &= \overline {\span} \{u_j\}_{j \in \bbN} \\
  &=  \overline {\span} (\{u_j\}^m_{j=1} \cup \{e_i : i \notin \bigcup
  \limits_{j=1}^m S_j\})\\
  &= \bigcap \limits^m_{j=1} \overline {\span} (\{u_j\} \cup \{e_i : i \notin
  S_j\})\\
  &= \bigcap \limits^m_{j=1} \bigcap \limits_{k \in S_{j} \backslash
  \{\varphi(j)\}} f^{-1}_k (0)
\end{align*}
and $(ii)$ is proven.

$(ii) \Longrightarrow (i):$
Suppose that there exist $k>1$ and $1 \leq i_1 < i_2 < \cdots < i_k \leq n$
such that $\card (\bigcup \limits^k_{\nu=1} F_{i_{\nu}})\le k$.  
Since
functionals $\{f_{i_{\nu}}\}^k_{\nu=1}$ are linearly independent the
 matrix of their coefficients 
$( f_{i_{\nu},j})^k_{\nu=1, j\in \cup^k_{\nu=1} F_{i_{\nu}}}$ 
has rank $k$. Thus $\card (\bigcup \limits^k_{\nu=1} F_{i_{\nu}})\ge k$.

If $\card (\bigcup \limits^k_{\nu=1} F_{i_{\nu}})= k$ then, again by the
linear  independence of the 
functionals $\{f_{i_{\nu}}\}^k_{\nu=1}$,  there exists
an isomorphism which transforms the $k \times k$ matrix of the 
coefficients of
$\{f_{i_{\nu}}\}^k_{\nu=1}$ into an identity matrix 
(of size $k \times k$), 
that is there exist functionals $\{g_{i_{\nu}}\}^k_{\nu=1}$ so that
for all $\nu=1, \cdots, k$ we have $\card(\supp
g_{i_{\nu}}) = 1$  and
$$
\bigcap \limits^k_{\nu=1} g_{i_{\nu}}^{-1} (0) = \bigcap \limits^k_{\nu=1}
f^{-1}_{i_{\nu}} (0).
$$

Thus without reducing the generality, we will assume that sets
$\{F_j\}^n_{j=1}$ satisfy the following property:
\begin{equation}
\lb{supports}
For \ each \ subset \ S \subset \{1, \cdots, n\} \ with \ \card(S)>1, \ \card(\bigcup
\limits_{i \in S} F_i) > \card (S)
\end{equation}
After reordering of $\{F_j \}^n_{j=1}$, if necessary, let $m$, 
$0 \leq m \leq
n$, be such that
$$
\card (F_j) = 2 \ for \ all \ j \leq m
$$
$$
\card(F_j) =1 \ for \ all \ j>m.
$$

Then $(\bigcup \limits^m_{j=1} F_j) \cap (\bigcup \limits^n_{j=m+1} F_j)=
\emptyset$ and
\begin{equation}
\lb{singles}
\bigcap \limits^n_{j=m+1} f^{-1}_j (0) = \overline {\span} \{e_k: k \notin
\bigcup \limits^n_{j=m+1} F_j\}.
\end{equation}

We introduce a relation $\sim$ on a set of indices $\{1, \cdots, m \}$ as
follows: $i \sim j$ if there exists $t \leq m$ and $k_1, \cdots, k_t \in \{1,
\cdots, m\}$ such that $i=k_1, j=k_t$ and $F_{k_{s}} \cap F_{k_{s+1}} \neq
\emptyset$ for all $1 \leq s < t$.

Clearly $\sim$ is an equivalence relation.  Let $S \subset \{1,2, \cdots, m\}$ be
a class of equivalence of $\sim$, and set 
$A_S = \bigcup \limits_{j \in S} F_j$.

Notice that
\begin{equation}
\lb{cardinality}
\card (A_S) =\card(S) +1
\end{equation}

Indeed, since $S$ is a class of equivalence of $\sim$, there exists a
bijection $\sigma: \{1, \cdots, \card S\} \to S$ such that $F_{\sigma(i)}
\cap F_{\sigma(i +1)} \neq \emptyset$  for each $i< \card(S)$.  Thus for $k
< \card(S)$ we have
\begin{align*}
\card(\bigcup \limits^{k+1}_{i=1} F_{\sigma(i)}) &\leq \card(\bigcup
\limits^{k}_{i=1}F_{\sigma(i)}) +1 \\
&\leq \card (F_{\sigma(1)}) +k =2+k
\end{align*}

Thus
$$
\card(A_S)= \card(\bigcup \limits^{\card(S)}_{i=1}F_{\sigma(i)}) \leq
\card(S) +1.
$$

On the other hand, by \eqref{supports}, if $\card(S)>1$
$$
\card (A_S) > \card(S).
$$

Thus $\card(A_S) = \card(S) +1$ if $\card(S)>1$.

If $\card(S)=1$, then $S=\{s\}$ for some $s \in \{1, \dots,m\}$ and $A_S =
F_s$, so $\card(A_S) = \card(F_s) =2$, and \eqref{cardinality} holds.

Now we are ready to show that for each class of abstraction of $\sim$, $S$,
there exists a vector $u_S \in X$ with $\supp u_S \subset A_S$ such that
\begin{equation}
\lb{classofabstraction}
\bigcap \limits_{j \in S} f^{-1}_j (0) = \overline{\span}(u_S \cup \{e_k: k
\notin A_S\})
\end{equation}

To see this notice that, clearly, if $k \notin A_S$ then $e_k \in
\bigcap \limits_{j \in S} f^{-1}_j (0)$.  Next let us consider $u$ with 
$\supp (u) \subset A_S$, say
$$
u = \sum \limits_{i \in A_S} u_i e_i\ .
$$

Then $u \in \bigcap \limits_{j \in S} f^{-1}_j (0)$ if and only if
$\{u_i\}_{i \in A_{S}}$ is a solution of a system of linear homogenous
equations $(f_j(u)= 0, \ j\in S)$ where the number of equations is $\card(S)$
and the number of variables is $\card(A_S)$.  But by \eqref{cardinality},
$\card (A_S) = \card(S)+1$ and the functionals $\{f_j\}_{j \in S}$ are
linearly independent, so this system has exactly one solution, which we
will denote $u_S$ and \eqref{classofabstraction} holds.  Combining
\eqref{classofabstraction} with \eqref{singles} we obtain
$$
\bigcap \limits^{n}_{j=1} f_j^{-1}(0)= \overline{\span} \{ \{e_k: k \notin \bigcup
\limits^n_{j=1} F_j\} \cup \{u_S: S \subset \{1, \cdots, m\}, S \ {\text 
{class 
of  abstraction   of} }\ \sim\}
$$
so $(i)$ is proven.
\end{pf}

\subsection{Characterizations of $\ell_p$ through 1-complemented 
subspaces} \lb{subseccharlp}
We finish this section with  results that generalize Ando's
characterization of $L_p$-spaces among Banach lattices
(Theorem~\ref{andolattice}) which were obtained by Calvert
and Fitzpatrick in a series of papers
\cite{CF86,CF87,CF87a,CF88,Cal88,CF89,FC91,FC91a,FC91b}.
The first result shows that Theorem~\ref{BP} characterizes $\ell_p$ and
$c_0$:

\begin{th}
\cite[Theorem~1]{CF86}
\lb{CF1}
Let $\{e_i \}^{\infty}_{i=1}$ be a Schauder basis for a Banach lattice $X$
with $e_i \wedge e_j =0$ if $i \neq j$.  Suppose each hyperplane $F$ which
contains all but two of the basis vectors (i.e. $F=f^{-1}(0)$ for some
functional $f$ with at most two coordinates different from $0$) and which
is a sublattice is the range of a contractive projection.  Then $X=\ell_p$ 
$(1 \leq p < \infty)$  or $X=c_0$.
\end{th}

Next they obtained a generalization of Ando's theorem about
$1$-complementability of two dimensional sublattices
\cite{ando69}, \cite[Theorem~16.4]{Lacey}:

\begin{th}
\cite[Theorem~4.2]{CF87} (cf. also \cite[Theorem~2]{CF87a})
\lb{CF2}
Let $X$ be a Banach lattice and $A=\{e_i\}_{i \in I}$ be a set of elements
of $X$ with $e_i \wedge e_j =0$ for $i \neq j $ and such that $\span
\{e_i\}_{i \in I}$ is dense in $X$.  Let $i_0 \in I$ be any fixed element
of $I$.  Suppose that any two dimensional sublattice of $X$ which contains
some $e_i, \ i \in I \backslash \{i_0\}$, is the range of a contractive
projection.  Then $X$ is linearly isometric and lattice isomorphic to
$\ell_p(I) \ (p \in [1, \infty])$ or $c_0(I)$.
\end{th}

Calvert and Fitzpatrick also show that in the assumptions of the above
theorem one cannot exclude two elements of the index set $I$, i.e. there
exists a Banach lattice $X \neq \ell_p (I),\  c_0(I)$ and $i_0, \ i_1 \in I,
\ i_0 \neq i_1$, such that every $2$-dimensional sublattice of $X$ containing
some $e_i,\  i \in I \backslash \{i_0, i_1\}$, is $1$-complemented
\cite[Example~4.]{CF87}.

Further, Calvert and Fitzpatrick obtained a result analogous to
Theorem~\ref{CF2} where the assumption about $X$ being a lattice is
replaced by the assumption of existence of enough smooth points in $X$.

\begin{th}
\cite[Theorem~C]{CF88}
\lb{CF3}
Let $X$ be a real Banach space with $\dim (X) \geq 3$.  Let $\{e_i\}_{i \in
I}$ be a linearly independent set of smooth points in $X$ with $\overline
{\span} \{e_i\}_{i \in I}=X$.  Suppose that every two-dimensional subspace
of $X$ intersecting $\{e_i\}_{i \in I}$ is the range of a nonexpansive
projection.  Then $X$ is isometrically isomorphic to $\ell_p (I),\  1 \leq p
\leq \infty$ or $c_0(I)$.
\end{th}

Finally in \cite{CF89} Calvert and Fitzpatrick studied extensions
of Theorems~\ref{CF2} and \ref{CF3} for spaces which are not lattices
and without the assumption of the existence of a linearly dense set
of smooth points in $X$ (like in Theorem~\ref{CF3}).
They showed that there exist 
$3$-dimensional real Banach spaces $X$, $X\neq \ell_p(3)$,
 which have two linearly independent
elements $e_1, \ e_2$ such that every $2$-dimensional subspace of $X$
containing either $e_1$ or $e_2$ is $1$-complemented in $X$. Moreover
they  gave a complete characterization of
$3$-dimensional real Banach spaces $X$ with the above property.
In particular, as a corollary of \cite{CF89} one obtains the following:

\begin{cor}
\lb{CF4}
Let $X$ be a real symmetric Banach space with $\dim (X) = 3$ such that
$X$ contains  two linearly independent
elements $e_1, \ e_2$ such that every $2$-dimensional subspace of $X$
containing either $e_1$ or $e_2$ is $1$-complemented in $X$. 
  Then $X$ is isometrically isomorphic to $\ell_p (3),\  1 \leq p
\leq \infty$.
\end{cor}

We do not know whether this corollary can be extended for symmetric spaces
$X$ with $\dim X >3.$

\begin{rem}
It seems that the  property like in Theorem~\ref{lp}(2) i.e.
\begin{equation} \lb{P}
\begin{split}
&\text{\it An infinite dimensional subspace $Y \subset X$ is 
1-complemented in $X$}\\
&\text{\it if and only if $Y$  is isometrically isomorphic to $X$.}
\end{split}
\end{equation}
does not characterize $\ell_p$.  Indeed Gowers and Maurey \cite{GM93}
constructed spaces which have very few complemented
subspaces, and we believe that these spaces can be slightly modified so 
that they would satisfy property \eqref{P} vacuously 
(cf. also Section~\ref{nonexistence} and \ref{seciso}).
\end{rem}

\subsection{Bibliographical remarks}        
 Other paper containing results about contractive projections in $\ell_p$ 
include \cite{Kin84,BP89,B90}.
Contractive projections in $\ell_p$ for $0<p<1$ were
 studied in \cite{Popa79,Popa81}.

\section{Sequence Banach spaces}
\lb{sequence}

\subsection{General nonexistence results} \lb{nonexistence}
The situation about the existence of contractive projections in Banach
spaces with bases is somewhat ambivalent.  That is, one may say that there
is an abundance of contractive projections since every infinite dimensional
Banach space with a basis has a conditional basis, which can be renormed
so that the projections associated with this conditional basis are all
contractive \cite[p.250]{Singer1}. 
Moreover Lindenstrauss proved the following property of nonseparable
reflexive Banach spaces illustrating the richness of 
1-complemented subspaces:

\begin{th}\cite{L64b}
Let $X$ be any reflexive Banach space.  For any separable subspace
$Y\subset X$ there exists a separable subspace $Z\subset X$ such that
$Y\subset Z$ and $Z$ is 
$1$-complemented in $X$.
\end{th}

On the other hand, we mentioned above
a striking example of Lindenstrauss illustrating the claim that norm
one projections are very rare, see Theorem~\ref{L64}.
Also
 Bohnenblust showed that there are
subspaces $S$ of $\ell_p^n$ which do not have any nontrivial 1-complemented
subspaces (Theorem~\ref{Bh2} above).
Moreover
Bosznay and Garay
\cite{BG} showed:

\begin{th}
\lb{BG}
Let $X$ be a finite dimensional real Banach space.  Then for any $\e
> 0$ there exists a norm $||| \cdot |||$on $X$ such that
$$
(1-\e) \| x \|_X \leq ||| x ||| \leq (1 + \e) \| x \|_X
$$
and such that $(X,||| \cdot |||)$ does not have any nontrivial
$1$-complemented subspaces.
\end{th}

Here, we say that a $1$-complemented subspace $F$ of $X$ is nontrivial if
$F \subsetneq X$ and $\dim F > 1$ (clearly, by Hahn-Banach theorem all
$1$-dimensional subspaces are $1$-complemented in all Banach spaces).

Thus, in the case of finite dimensional real spaces, the property of having
$1$-complemented subspaces is highly unstable.  It seems clear that there
should exist infinite dimensional Banach spaces with no nontrivial
$1$-complemented subspaces (we believe that an appropriate modification of a
hereditarily indecomposable space constructed in \cite{GM93} would do the
job).  However the following questions are open:

\begin{question} \lb{renorm}
Let $(X, \| \cdot \|)$ be a Banach space.  Does there exist an equivalent
renorming $(X, ||| \cdot |||)$ of $X$ so that $(X, ||| \cdot |||)$ has no
nontrivial $1$-complemented subspaces?
\end{question}

Theorem~\ref{BG} says that when $X$ is finite dimensional then the answer
to this question is positive in the very strong sense that one may
require the new norm $|||\cdot |||$ to be arbitrarily close to the original
norm $\| \cdot \|$. 

\begin{question} \lb{subspace}
Let $X$ be a Banach space.  Does there exist a subspace $Y \subset X$ so
that $Y$ does not have nontrivial 1-complemented subspaces?
\end{question}

As mentioned above, by a result of Bohnenblust,
 Question~\ref{subspace} has a positive answer if
$X=\ell_p^n, \ p\in(1,\infty)\setminus \bbZ, \ n\ge 7,$ and therefore 
also if $X=L_p$ or $X=\ell_p$. We suspect that the answer is positive
for all Banach spaces which are not isometric to a Hilbert space.

\subsection{More negative results concerning the inheritance of the isomorphic
structure}
Here we will concentrate on the search of characterizations of
$1$-complemented subspaces and norm one projections in classical types of
spaces (with their usual norms) in the spirit of Theorem~\ref{lp}.
However we have to start from the following two negative results:
\begin{th}
\cite{P71}
\lb{P71}
Let $X$ be any finite dimensional Banach space $(\dim X=n)$. Then
\begin{itemize}
\item[(a)]
there exists a Banach space $Y$ with $\dim Y=n^2$ such that $Y$ has a basis
with basis constant $\leq 1 + n^{- 1/2}$ and $X$ is a
$1$-complemented in $Y$.
\item[(b)]
there exists a Banach space $Y$ such that $Y$ has a monotone basis (i.e.
basis constant equals $1$, cf. Definition~\ref{defbasis}) and $X$ is 
$1$-complemented in $Y$ (here $\dim Y$
could be infinite).
\end{itemize}
\end{th}

\begin{th}
\cite{L72}, \cite[Theorem~3.b.1]{LT1}
\lb{L}
Every space $X$ with a $1$-unconditional basis is $1$-complemented in some
symmetric space $X$.
\end{th}
These results are very negative because they show that $1$-complemented
subspaces do not have to inherit the isomorphic structure of the space.  In
particular, it follows from Theorem~\ref{L}, that there is no hope of
extending Proposition~\ref{isometryH} to general symmetric spaces
(Proposition~\ref{isometryH} does have an analogue in $\ell_p$,
Theorem~\ref{lp}(2)).

Moreover very shortly after Theorem~\ref{L}
it was established that in fact the isomorphic structure of
$X$ is not inherited by $1$-complemented subspaces even in the most
natural symmetric spaces i.e. in Orlicz spaces and Lorentz spaces
(cf. Definitions~\ref{defo} and \ref{defl}).  We have:

\begin{th}
\cite{L73,ACL73}
\lb{nonisomorphic}
There exists a class of Orlicz spaces $\ell_\varphi$ \cite{L73} and a class
of Lorentz spaces $\ell_{w,p}$ \cite{ACL73} such that if $X$ belongs to
either of these classes then $X$ has an infinite dimensional
$1$-complemented subspace $Y$ which is not isomorphic with $X$.
\end{th}

Thus it seems that there is no hope of giving any sort of characterization of
$1$-complemented subspaces in terms of isomorphisms.
However the analogue of Theorem~\ref{condexp} still holds:

\begin{th}
\lb{block}
Let $X$ be a symmetric sequence space with basis $\{e_i\}_{i \in I} \ (I
\subseteq \Bbb N)$.  Let $\{f_j\}_{j \in J}$ be a block basis with constant
coefficients of any permutation of $\{e_i\}_{i \in I}$, i.e. there exists a
permutation $\sigma$ of $I$, an increasing sequence $\{p_j\}\subseteq
\Bbb N$ and scalars $\{\theta_i\}_i$ with $|\theta_i|=1 $ for all $i\in I$
 such that
$$
f_j = \sum \limits^{p_{j+1}}_{k=p_{j}+1} \theta_ke_{\sigma(k)}.
$$
Then $Y= \overline {\span} \{f_j\}_{j \in J} \subset X$ is a
$1$-complemented in $X$ and the averaging projection $P: X \ONTO Y$ is a
norm-one projection.
\end{th}

In fact all examples presented in proofs of Theorems~\ref{L} and
\ref{nonisomorphic} were of the form described in Theorem~\ref{block}.  We
postulate that Theorem~\ref{block} is the right form of the description of
$1$-complemented subspaces, i.e. that Conjecture~\ref{condexpcon} should be
valid also in Banach spaces with bases.

\begin{con}
\lb{blockcon}
Let $X$ be a strictly monotone Banach space with a $1$-unconditional basis
(sufficiently different from a Hilbert space).  Then every norm one
projection in $X$ is a weighted conditional expectation operator and every
$1$-complemented subspace of $X$ is spanned by mutually disjoint elements.
\end{con}

Below we present some results supporting this conjecture (see
Theorem~\ref{BR1} and Corollary~\ref{BR1cor}), but first we return to the
chronological order of discoveries.

\subsection{Some spaces whose 1-complemented subspaces do inherit the basis}
\lb{ecp}
The next developments in the study of $1$-complemented subspaces
of sequence space dealt with isometric versions of Banach's 
Problems~\ref{basis}
 and
\ref{comp}. In this section we will discuss the isometric analogue of
Problems~\ref{basis} i.e.:
\begin{problem}
\lb{monotone}
Does every $1$-complemented subspace of a space with a monotone basis 
have a monotone basis? (cf. Definition~\ref{defbasis})
\end{problem}

As we mentioned above, Pe\l czy\'nski \cite{P71} (see Theorem~\ref{P71}(b)) showed
that the answer to Problem~\ref{monotone} is negative.  However Dor
\cite{Dor73} proved that the answer is yes if we consider only finite
dimensional spaces:

\begin{th}
\lb{Dor73}
Let $X$ be a finite dimensional Banach space with a monotone basis, and let
$Y$ be $1$-complemented in $X$.  Then $Y$ has a monotone basis.
\end{th}

Chronologically the next development related to describing  bases
in  1-complemented subspaces is due to Kinnunen \cite{Kin84} which we
described in Section~\ref{ortho}. 

Next development concerning Problem~\ref{monotone} is due to Rosenthal
\cite{R84} who proved:

\begin{th}
\lb{R84}
Let $Y$ be a reflexive Banach space which is isometric to a contractively
complemented subspace of a Banach space $X$ with reverse monotone basis.
Then $Y$ has a reverse monotone basis.
\end{th}

Here, a basis $(b_j)$ for a Banach space $X$ is said to be {\it reverse
monotone} if $\| I-P_j \| =1$ for all $j$, where
$$
P_jb = \sum \limits^j_{n=1} c_n b_n \ \ {\text{ if }} \ \ b = \sum^\infty_{n=1} c_n
b_n.
$$

Clearly if $X$ is finite dimensional then $X$ has a reverse monotone basis
if and only if $X$ has a monotone basis (simply reverse the order of the
monotone basis $(b_n)_{n=1}^m$ to $(b_n)^1_{n=m})$.

Thus Theorem~\ref{R84} generalizes Theorem~\ref{Dor73} and also it
illustrates that in infinite dimensions the concept of a reverse monotone
basis is very different from a monotone basis (cf. Theorem~\ref{P71}(b)).

In \cite{R84} Rosenthal also studied the following concept:

\begin{defn}
A Banach space $X$ has enough contractive projections $(ECP)$ provided
every $1$-complemented non-zero subspace $Y$ of $X$ contains a
contractively complemented subspace $Z$ of codimension one in $Y$.

Property $(ECP)$ is clearly inherited by $1$-complemented subspaces.
\end{defn}

This property seems to be modeled on the definition of the reverse
monotone basis $(b_n)_{n \in \Bbb N}$ where we require that every subspace
of form $Y_M = \overline{\span} \{b_n\}_{n \geq M}\  (M \in \Bbb N)$ has a
$1$-complemented subspace $Y_{M+1}=\overline{\span} \{b_n\}_{n \geq M +1}$
of codimension one in $Y_M$.

 Rosenthal proved:
\begin{th}
\cite[Theorem~1.8]{R84}
\lb{ECP}
Every Banach space $X$ with reverse monotone basis has enough contractive
projections $(ECP)$.
\end{th}

The proof of Theorem ~\ref{ECP} is not trivial. It uses the theory of
numerical ranges (see \cite {BD1,BD2}) and Proposition~\ref{flinn}.

On the other hand Rosenthal suggests that there may exist reflexive
separable spaces with $ECP$ but with no basis or even without
finite-dimensional decomposition $(FDD)$. This question is still open.

Rosenthal  proved the following characterization of property $(ECP)$:

\begin{th}
\cite[Theorem~2.1]{R84}
A reflexive Banach space $X$ has $(ECP)$ if and only if $X$ has a reverse
monotone transfinite basis.
\end{th}

Here transfinite basis is the concept that generalizes bases by dropping the
assumption of countability (it was introduced by Bessaga \cite{Be67}, see
also \cite{Dor71}).  It is defined as follows:

\begin{defn}
\cite[Definition~17.7]{Singer2}
Let $\eta > 0$ be an ordinal.  A transfinite sequence of elements
$(b_\alpha)_{\alpha < \eta}$ in a Banach space $X$ is called a {\it
transfinite basis} (of length $\eta$) of $X$ if for every $x \in X$ there
exists a unique transfinite sequence of scalars $(x_\alpha)_{\alpha <
\eta}$ such that $\sum \limits_{\alpha < \eta} x_\alpha b_\alpha$ converges
to $x$.
\end{defn}

For most recent developments related to further generalizations of 
Problem~\ref{monotone} see Section~\ref{ap}.

\subsection{Preservation of the 1-unconditional  basis in the complex case}
Let us now concentrate on the isometric version of Problem~\ref{comp} i.e.:
\begin{problem}
\lb{1-unc}
Does every $1$-complemented subspace of a space with a $1$-unconditional basis
have  a $1$-unconditional basis?
\end{problem}

Curiously the situation is different depending whether we consider Banach
spaces over $\Bbb C$ or over $\Bbb R$.

The first result concerns complex spaces and was proven implicitly by Kalton
and Wood \cite {KW} and explicitly explained in \cite {F84,R83}.  We have

\begin{th} \lb{KW}
Let $X$ be a complex Banach space with 1-unconditional basis and let $Y$ be
a 1-complemented subspace of $X$.  Then $Y$ has a 1-unconditional basis.
\end{th}

The analogous result is false in real Banach spaces \cite{L79,BFL84}. Thus
the answer to Problem~\ref{1-unc} is negative in the real case.

The proof of Theorem~\ref{KW} is based on the theory of numerical ranges
in Banach spaces (see \cite{BD1,BD2}).  We outline here the main elements
of the proof.

\begin{defn}
\cite{BD1}
Let $X$ be a Banach space and $T$ be a linear operator on $X$. We say that
operator $T$ is {\it hermitian}  if numerical range of $T$ is contained in
$\Bbb R$.  i.e.
$$
\overline{\conv}\{\left< f, T x \right>  : x \in X, f \in J(x)\} \subseteq \Bbb R
$$
\end{defn}

\begin{defn} \cite{KW}
An element $x_0 \in X$ is called {\it hermitian} if there exists a
hermitian rank-1 projection from $X$ onto $\span\{x_0\}$, i.e. equivalently,
if for all $ x \in X, f \in J(x), f_0 \in J (x_0)$ we have
$$
\left< f,  x_0 \right> \cdot \left<f_0, x\right> \in \Bbb R.
$$

The set of all hermitian elements of $X$ is denoted $H(X)$

Let $\{H_\lambda : \lambda \in \Lambda\}$ be the collection of maximal
linear subspaces of $H(X)$.  Then $H_\lambda, \ \lambda \in \Lambda$ are
called {\it Hilbert components} of $X$.

A Hilbert component $H_\lambda$ is called {\it nontrivial} if
$\dim H_\lambda > 1 $
\end{defn}

Kalton and Wood showed that Hilbert components are well-defined and they
obtained the following characterization of hermitian elements in $X$:

\begin{th} \cite [Theorem ~6.5]{KW} \lb{KW2}
Let $X$ be a complex Banach space with a normalized 1-unconditional basis
$\{e_i\}_{i\in I}$.
Then $x_0 \in X$ is hermitian if and only if
\begin{itemize}
\item [(i)] $\Vert y \Vert_X = \Vert y \Vert_2$ for all $y \in X$ with $\supp
y \subset \supp x_0$, and
\medskip
\item [(ii)]  for all $y,z \in X$ with $\supp y \cup  \supp z \subset \supp
x_0$ and for all $v \in X$ with $\supp v \cap \supp x_0 = \emptyset$ if
$\Vert y \Vert_X = \Vert z \Vert_X$ then $\Vert y + v \Vert_X =
\Vert z + v \Vert_X$.
\end{itemize}
\end{th}

In (i) above when we say $\|y\|_2 $ we mean the $\ell_2-$norm of $y$ with respect
to the given 1-unconditional basis of $X$, i.e. if $y=\sum_{i\in I} y_ie_i$
then $\|y\|_2=(\sum_{i\in I}|y_i|^2)^{1/2}.$

In particular, it is easy to see that
\begin{cor}
\lb{KW2cor}
If $X$ is a complex Banach space with a normalized 1-unconditional basis
$\{e_i\}_{i \in I}$ and $x_0$ is hermitian in $X$ then
$$
S_{x_0} = \overline{\span} \{e_i : i \in \supp x_0\} \subset X
$$
is a 1-complemented subspace of $X$ and $S_{x_0}$ is isometric to a Hilbert
space.
\end{cor}

Clearly every element $\{e_i\}_{i \in I}$ of the 1-unconditional basis is
hermitian, and Kalton and Wood proved that if $P$ is a projection of norm
one then also the elements $\{P{e_i}\}_{i \in I}$ are hermitian in $\calR(P)$ (see
\cite[Lemma~3]{F84}).  This fact, together with Theorem~\ref{KW2} leads to
Theorem~\ref{KW}.

A refinement of this technique allowed the author of this survey to obtain
a strengthening of Theorem~\ref{KW}.  Namely we get the geometric description of the
1-unconditional basis of the 1-complemented subspace:

\begin{th}
\cite[Corollary~3.2]{complex}\lb{BR1}
Suppose that $X$ is a complex strictly monotone Banach space with a
1-unconditional basis $\{e_i\}_{i}$ and let $P$ be a projection of norm one
in $X$.  Suppose that $Y = P(X) \subset X$ has no nontrivial Hilbert
components.  Then there exist disjointly supported elements $\{y_i\}^{\
m}_{i=1} \subset Y \ (m = \dim Y \leq \infty)$ which span $Y$.  Moreover,
for all $x \in X$
$$
Px = \sum^m_{i=1}y_i^*(x)y_i
$$
where $\{y_i^* \}^{\ m}_{i=1} \subset X^*$  satisfy $\Vert y_i^*\Vert =
\Vert y_i\Vert_X = y_i^* (y_i) = 1$ for all $i = 1,\ldots,m$.
\end{th}

The statement in Theorem~\ref{BR1} exactly parallels and extends the
characterization of 1-complemented subspaces of $\ell_p$ given in
Theorem~\ref{lp}(3). 

\begin{rem}
The assumption of $X$ being strictly monotone cannot be omitted (Blatter
and Cheney \cite{BCh} showed examples of 1-complemented hyperplanes in
$\ell^3_\infty$ that are not spanned by disjointly supported vectors).
Also the assumption that $Y$ does not have non-trivial Hilbert components
cannot be dropped (see examples in \cite{complex}); the property of
not having non-trivial Hilbert components is not necessarily inherited by
1-complemented subspaces.
\end{rem} 

As noted above (Theorem~\ref{nonisomorphic}) one
cannot hope to give any sort of isomorphic description of 1-complemented
subspaces of $X$ and we believe that the geometric characterization of the
position of $Y$ in $X$ as the span of a family of disjointly supported
vectors is best possible.  In fact the examples of 1-complemented subspaces
$Y$ nonisomorphic to the whole space $X$ which were described in
\cite{ACL73, L73} all satisfied the conclusion of Theorem~\ref{BR1} (in
both complex and real case), even more, they in fact were spanned by
disjoint elements with finite supports and constant coefficients i.e. were
of the form
$$
Y = \overline{\text{span}} \{v_j\}_{j \in \Bbb N}
$$
where $v_j = \sum_{\nu\in S_j}e_\nu,\ S_j \cap S_k = \emptyset$ whenever $j
\not= k$, so $Y$ satisfies the conditions of Theorem~\ref{LT} and the norm one
projection onto $Y$ is simply a conditional expectation operator.

Thus it appears that it is not that the class of contractive projections is
richer in spaces other than $\ell_p$, but rather that the isomorphic
structure of subspaces of $X$ is much more varied.  In fact,
Theorem~\ref{BR1} can be rephrased as follows:

\begin{cor} \lb{BR1cor}
Suppose that $X$ is a complex Banach space with 1-unconditional basis and
$X$ does not contain a 1-complemented copy of a 2-dimensional Hilbert space
$\ell^2_2$.  Then every norm one projection in $X$ is a weighted
conditional expectation operator.
\end{cor}

\begin{pf}
Theorem~\ref{BR1} and Corollary~\ref{KW2cor}.
\end{pf}

\begin{rem}
Recall that Calvert and Fitzpatrick showed that if all weighted conditional
expectation operators are contractive projections in $X$, then $X$ has to be isometric
to $\ell_p$ or $c_0$ (see Section~\ref{lpsection}).
\end{rem}

\subsection{The real case}
One can immediately see that when $X$ is a real Banach space then all
operators and all elements are hermitian, so above methods cannot be
applied.  And, in fact, Theorem~\ref{KW} fails in the real case.

As mentioned above, Lewis \cite{L79} and Benyamini, Flinn, Lewis
\cite{BFL84} showed examples of real Banach spaces without 1-unconditional
bases which are 1-complemented in a real space with a 1-unconditional basis:

\begin{prop}
\lb{example}\cite[Propositions~1 and 2]{BFL84}
The space $E_n = \{(x_i)^{\ n}_{i=1} \in \ell_\infty^n : \sum^n_{i=1}x_i =
0\}$ is 1-complemented in a space with 1-unconditional basis and if $n \geq
5$ then $E_n$ does not admit a 1-unconditional basis.
\end{prop}

The proof of Proposition~\ref{example} is based on the following
interesting observation:

\begin{prop} (stated as ``well-known'' in \cite{BFL84})
\lb{obser}
A real Banach space $X$ embeds as a norm-one complemented subspace of some
Banach space $Z$ with 1-unconditional basis if and only if the identity $I$
on $X$ can be written as
$$
I = \sum_i f_i \otimes x_i
$$
\noindent
for some $x_i \in X, f_i \in X^*$ with
$$
\Vert \sum_i \varepsilon_i f_i \otimes x_i \Vert \leq 1
$$
\noindent
for all $\varepsilon_i = \pm 1$.
\end{prop}

In Proposition~\ref{example} spaces $E_n$ are not uniformly convex, but
uniformly convex examples are also possible.

It is easy to see that spaces $E_n$ do have a 2-unconditional basis.  The
following version of Problem~\ref{comp} is open:

\begin{problem}\lb{1comp}
Let $X$ be a real Banach space with a 1-unconditional basis and let $Y$ be
1-complemented in $X$.
Does $Y$ have to have an unconditional basis?
\end{problem}

After the negative examples of Lewis and Benyamini, Flinn, Lewis there were
few attempts to characterize 1-complemented subspaces of special real
Banach spaces.  First development is due to Rosenthal \cite {R86} who
considered spaces which are isometric to the direct sum of Hilbert spaces
of dimension at least two via a one-unconditional basis according to the
following definition:

\begin{defn}\cite{R86}
Let $\Gamma$ be a nonempty set and $(X_\alpha)_{\alpha \in \Gamma}$ be a
family of nonzero Banach spaces.  A Banach space $X$ is said to be a {\it
functional unconditional sum of the $X_\alpha$'s} if there exists a
normalized 1-unconditional basis $\underline b = (b_\alpha)_{\alpha \in
\Gamma}$ for some Banach space $B$ so that $X$ is (linearly isometric
to) ${\disp(\sum_\Gamma \oplus X_\alpha)_{\underline b}}\ $,
i.e.
$$
X = \{x = (x_\alpha)_{\alpha \in \Gamma} \in \prod_{\alpha \in \Gamma}
X_\alpha : \sum_{\alpha \in \Gamma} \Vert x_\alpha\Vert_{X_\alpha}b_\alpha
\in B\}
$$
and we define the norm on $X$ by:
$$\Vert x\Vert_X = \Vert \sum\limits_{\alpha \in
\Gamma}\Vert x_\alpha b_\alpha \Vert _{B}.$$

In case $X$ is real and each $X_\alpha$ is a real Hilbert space of
dimension at least two, we call $X$ a {\it Functional Hilbertian Sum}
(FHS).
\end{defn}

It is clear that each complex Banach space $X$ with a 1-unconditional basis
$(e_i)_{i \in I}$
is isometrically isomorphic over reals to $X_{\Bbb R}(\ell_2^2)$  
$X_{\Bbb R} = \{\sum_{i \in I} a_ie_i \in X :$ where  
 $(a_i)_{i \in I} \subset \Bbb R\}$  i.e.
$X$ can be considered over reals as an FHS where $B = X_{\Bbb R},\
\underline b = (b_i)_{i \in I}: = (e_i)_{i \in I}$ and each $(X_i)_{i\in
I}$ is a two-dimensional Hilbert space $\ell^2_2$.  This similarity to the
complex case allows to transfer some of the techniques related to hermitian
operators and hermitian elements.  We will use the following notion:

\begin{defn}
Let $B$ be a complex or real Banach space.  We say that operator $T:B
\lra B$ is {\it skew-Hermitian} if $\Re f(Tb) = 0$ for all $b \in B$ 
and $f \in J(b)\subset B^*$.
\end{defn}

Using skew-Hermitian operators, Rosenthal gave a sufficient condition for a
1-comp\-le\-ment\-ed subspace of a real Banach space with 1-unconditional basis
to have a 1-unconditional basis:

\begin{th}
\cite[Theorem~3.15]{R86} \lb{R1}
Let $X$ be a real Banach space with a 1-unconditional basis and $Y$ be a
1-complemented subspace of $X$.  Suppose for all $y \in Y$ there is a
skew-Hermitian operator $T$ on $X$ such that $Ty \in Y$ and $T^2 y = -y$.
Then $Y$ is FHS, so $Y$ has a 1-unconditional basis.
\end{th}

As a corollary Rosenthal obtains Theorem~\ref{KW}.  The method of Rosenthal
is very similar in spirit to the method of Kalton and Wood \cite{KW}, but
Rosenthal's analysis of FHS spaces also involves deep considerations 
of Lie
algebras of Banach space $X$ (the Lie algebra of $X$ consists of
skew-Hermitian operators on $X$).

Next Rosenthal considers orthogonal projections on $X$ which are defined as
follows:

\begin{defn}
Let $X$ be a real or complex Banach space, and $Y, Z$ be subspaces of $X$.
$Z$ is said to be an {\it orthogonal complement of $Y$} if $Y + Z = X$ and
if for all $y \in Y, z \in Z$ and all scalars $\alpha, \beta$ with $\vert
\alpha \vert = \vert \beta \vert = 1,\ \Vert y + z \Vert = \Vert \alpha y +
\beta z \Vert$.

The projection $P$ with range $Y$ and kernel $Z$ is called {\it the
orthogonal projection onto $Y$}.  Note that $\Vert P \Vert = \|I - P\Vert = 1$. We say that {\it $Y$ is orthogonally complemented in $X$ if $Y$ has an
orthogonal complement}
\end{defn}

It is easy to see that any space with 1-unconditional basis is orthogonally
complemented in some FHS space.  Rosenthal proved the converse:

\begin{th}\cite[Theorem~3.18]{R86}\lb{R2}
Let $X$ be an FHS space and $Y$ be an orthogonally complemented subspace of
$X$.  Then $Y$ has a 1-unconditional basis.
\end{th}

Theorems~\ref{R1} and \ref{R2} are the most general partial answers to
Problem~\ref{1comp} known today.

It is natural to ask whether the structures analogous to but weaker 
than 1-unconditional bases are preserved by 1-complemented subspaces.
Results of this kind will be presented in Section~\ref{ap} below.

Next results on $1$-complemented subspaces of real Banach spaces with
1-unconditional bases 
are due to
the author of this survey.  In \cite{rocky} we considered an extension of
Theorem~\ref{BP} from $\ell_p$ to a larger class of Banach spaces.

\begin{defn}
(cf. \cite[Definition~1.d.3]{LT2})
We say that a Banach lattice $X$ is {\it one-p-convex} (resp. {\it
one-p-concave}) if for every $n \in \Bbb N$ and every choice of elements
$\{x_i\}_{i=1}^{n}$ in $X$ we have
$$
\| (\sum \limits^{n}_{i=1} |x_i|^p)^{\frac{1}{p}} \| \leq (\sum
\limits^{n}_{i=1} \| x_i \|^{p})^{\frac{1}{p}} \ \text{ if } 
\ 1 \leq p < \infty,
$$
or, respectively,
$$
\|(\sum \limits^{n}_{i=1}| x_i |^{q})^{\frac{1}{q}} \| 
\geq (\sum
\limits^n_{i=1} \|x_i \|^{q})^{\frac{1}{q}} \ \text{ if }  
\ 1 \leq q < \infty.
$$
\end{defn}

This is the usual notion of $p$-convexity (resp. $q$-concavity) but with
the additional requirement that the constant is equal to $1$.  Clearly
spaces $\ell_p$ are both one-p-convex and one-p-concave.

We have the following extension of Theorem~\ref{BP}.

\begin{th}
\cite[Theorem~1]{rocky}
\lb{BRrocky}
Let $X$ be a strictly monotone Banach space with a $1$-unconditional basis.
Suppose that
\begin{itemize}
\item[(a)]
$X$ is a one-p-convex, $2<p<\infty$, or
\item[(b)]
$X$ is a one-q-concave, $1<q<2$, and smooth at each basis vector.
\end{itemize}

Then any $1$-complemented subspace $F$ of codimension $n$ in $X$ contains
all but at most $2n$ basis vectors of $X$, i.e. there exist functionals
$f_n, \cdots, f_n \in X^*$ such that $\card (\bigcup \limits^n_{j=1} \supp
f_j) \leq 2n$ and $F= \bigcap \limits^n _{j=1} f_j^{-1}(0)$.
\end{th}

Notice that in a general Banach space $X$ it does not have to be true that
every  $1$-complemented subspace
$Y\subset X$ with $\codim Y<\infty$ is
an intersection of $1$-complemented hyperplanes in $X$ (as is the case
when $X=\ell_p$, by Theorem~\ref{BP}).  Indeed Bohnenblust's finite  
dimensional
examples from Theorem~\ref{Bh2} do not have $1$-complemented hyperplanes,
while every $1$-dimensional subspace is of finite codimension and, of
course, $1$-complemented (moreover these spaces are  
one-p-convex when $2<p<\infty)$,
see also Theorem~\ref{BG}.

We wish to note here that the condition in the conclusion of 
Theorem~\ref{BP} does not characterize $\ell_p$ among infinite dimensional
Banach spaces with 1-unconditional basis (cf. \cite[Theorem~5.1]{complex}).
However we   believe that the following conjecture is true:

\begin{con}
\mbox{  }

\begin{itemize}
\item[(a)]
Let $X$ be a finite dimensional Banach space such that every 1-dimensional
subspace of $X$ can be represented as an intersection of 
$1$-complemented hyperplanes. Then $X$ is isometric to $\ell_p^n$ for some
$p, 1\le p\le\infty, \ n=\dim X$.
\item[(b)]
There exists  a finite dimensional Banach space $X$, $\dim X=n$, which is
not   isometric to   $\ell_p^n$ for any
$p, 1\le p\le\infty, $ and
such that every 1-complemented
subspace $Y$ of $X$ with $\dim Y\ge 2$ can be 
represented as an intersection of 
$1$-complemented hyperplanes. 
\end{itemize}
\end{con}

As a corollary of Theorem~\ref{BRrocky} and Theorem~\ref{hyperplanes} we
obtain the extension of Theorem~\ref{BM}:

\begin{cor}
\cite[Corollary~3]{rocky}
Suppose that $X$ is a separable strictly monotone real function space on
$(\Omega, \mu)$ where $\mu$ is finite on the nonatomic part of $\Omega$.
Suppose that $X$ is either one-p-convex for some $2<p<\infty$, or
one-q-concave for some $1<q<2$ and smooth at $\chi_A$ for every atom $A$
of $\mu$.  Let the hyperplane $F=f^{-1}(0)$ be $1$-complemented in $X$.
Then there exist $\alpha, \be \in \Bbb R$ and $A,B$-atoms of $\mu$ so that
$f=\alpha \chi_A + \be \chi_B$.
\end{cor}

Theorem~\ref{BRrocky} does not guarantee that the $1$-complemented
subspaces of finite codimension in $X$ are spanned by disjointly supported
vectors, however we believe that this should be the case in ``most'' spaces.
In \cite{real} we proved that this is the case for most $1$-complemented
subspaces of finite codimension in Orlicz sequence spaces and
$1$-complemented hyperplanes in Lorentz sequence spaces.

We need a following notation:

\begin{defn}
We say that the Orlicz function $\varphi$ is {\it similar} to $t^p$ for
some $p \in[1, \infty)$ if there exist $C, t_0 > 0$ so that
$\varphi(t)=Ct^p$ for all $t<t_0$.

We say that $\varphi$ is {\it equivalent} to $t^p$ for some $p
\in[1,\infty)$ if there exist $C_1,C_2,t_0 > 0$ so that $C_1,t^p \leq
\varphi(t) \leq C_2 t^p$ for all $t<t_0$.
\end{defn}

We obtained the following:

\begin{th}
\cite[Theorem~7]{real}
\lb{orlicz}
Let $\varphi$ be an Orlicz function such that $\varphi$ is not similar to
$t^2$ and $\varphi(t)>0$ for all $t>0$.
Consider the Orlicz space $\ell_\varphi$ with either Luxemburg or Orlicz
norm.  Let $F \subset \ell_\varphi$ be a subspace of codimension $n$ with
$\dim F > 1$.  If $F$ contains at least one basis vector and $F$ is
$1$-complemented in $\ell_\varphi$ then $F$ can be represented as $F=
\bigcap \limits^n_{j=1} f_j^{-1}(0)$ where $\card (\supp f_j) \leq2$ for
all $j=1, \dots, n$, i.e. $F$ is spanned by disjointly supported vectors.

Moreover, if $\varphi$ is not equivalent to $t^p$ for any $p \in[1,
\infty)$, then $|f_{ji}|$ is either $1$ or $0$ for all $i,j$, i.e. $F$ is a
span of a block basis with constant coefficients.

If $\varphi$ is equivalent but not similar to $t^p$ for some $p \in[1,
\infty)$, then there exists $\g \geq 1$ such that
$\{|f_{ji}|: j=1, \dots, n, i \leq \dim\ell_{\varphi}\}\subset \{\g^m : m
\in {\Bbb Z}\} \cup \{0\}$.
\end{th}

In \cite{real} we were unable to eliminate the condition that $F$ has to
contain at least one basis vector of $\ell_\varphi$ for
Theorem~\ref{orlicz} to hold, even though we believe that this condition is
not necessary.  It follows from Theorem~\ref{BRrocky} that this condition
will be automatically satisfied if $\dim(\ell_\varphi)$ is large enough and
$\ell_\varphi$ is one-p-convex $2<p<\infty$ or one-q-concave $1<q<2$.
Very recently we also obtained a new result which eliminates this condition
(as well as the requirement that $F$ is of finite codimension) provided
that Orlicz function $\varphi$ is smooth enough \cite{preparation}.  The
other conditions in the assumptions of Theorem~\ref{orlicz} are all
necessary.

In \cite{real} we also considered characterizations of $1$-complemented
hyperplanes in Lorentz sequence spaces.  We obtained the following:

\begin{th}
\cite[Theorem~3]{real}
\lb{lorentz1}
Let $\ell_{w,p}$ be a real Lorentz sequence space with $1<p<\infty$ and
$w_2>0$.  Suppose that $Y=f^{-1}(0)$ is $1$-complemented in $\ell_{w,p}$
and $\card (\supp f) \geq n > 2$.  Then $p=2$ and $1=w_1=w_2=\cdots=w_n$.
\end{th}

\begin{th}
\cite[Corollary~6]{real}
\lb{lorentz2}
Let $\ell_{w,p}$ be a real Lorentz sequence space with $1<p<\infty$ and
$w_k >0$ for all $k$, i.e. $\ell_{w,p}$ is strictly monotone.  Suppose that
$Y=f^{-1}(0)$ is $1$-complemented in $\ell_{w,p}$ and $\card(\supp f)=2$
i.e. $f=f_i e_i^k + f_j e_j^k$ for some $i \neq j$.  Then $|f_i| = |f_j|$
or $\ell_{w,p}=\ell_p$ i.e. $w_k = 1$ for all $k$.
\end{th}

The main tool in proofs of Theorem~\ref{orlicz}, \ref{lorentz1} and
\ref{lorentz2} is Proposition~\ref{flinn}, which allows to transfer to the 
real space setting techniques analogous to hermitian elements in 
complex spaces (see \cite{KR}).

Finally we want to mention two results valid in complex and real Orlicz and
Lorentz sequence spaces, which characterize which subspaces are
$1$-complemented among the subspaces spanned by a family of mutually
disjoint vectors.

\begin{th}
\cite[Theorem~6.1]{complex}
\lb{blocko}
Let $\ell_\varphi$ be a real or complex Orlicz space and let $\{x_i\}_{i
\in I} \ \ (I \subseteq {\Bbb N},\  \card(I)>1)$ be mutually disjoint 
elements
in $\ell_\varphi$ with $\| x_i \|_\varphi = 1$ for all $i \in I$.  Suppose
that $X = \overline {\span} \{x_i\}_{i \in I} \subset \ell_\varphi$ is
$1$-compelemented in $\ell_\varphi$.  Then one of the three possibilities
holds:
\begin{itemize}
\item[(1)]
for each $i \in I,\ \card (\supp x_i) < \infty$ and $|x_{ij}| = |x_{ik}|$
for all $j, k \in \supp x_i$; or
\item[(2)]
there exists $p,\  1 \leq p \leq \infty$, such that $\varphi(t) = C t^p$ for
all $t \leq \sup \{\| x_i \|_\infty : i \in I\}$; or
\item[(3)]
there exists $p, \ 1 \leq p \leq \infty$, and constants $C_1, C_2,
\gamma \geq 0$ such that $C_2t^p \leq \varphi (t) \leq C_1
t^p$ for all $t \leq \sup \{\| x_i \|_\infty: i \in I\}$ and such that for
all $i \in I$ and $j \in \supp x_i$
$$
| x_{ij} | \in \{\gamma^k \cdot \|x_i \|_\infty: k \in \Bbb Z\}.
$$
\end{itemize}
\end{th}

\begin{th}
\cite[Theorem~6.3]{complex}
\lb{blockl}
Let $\ell_{w,p}$ with $1<p<\infty$, be a real or complex Lorentz sequence
space and let $\{x_i\}_{i \in I} (I \subseteq {\Bbb N}, \card (I)>1)$ be
mutually disjoint elements in $\ell_{w,p}$ such that $X = \overline {\span}
\{x_i\}_{i \in I}$ is $1$-complemented in $\ell_{w,p}$.  Suppose, moreover,
that $w_\nu \neq 0$ for all $\nu \leq S \overset{ \text{def} }{=} \sum
\limits_{i \in I} \card (\supp x_i)\  (\leq \infty)$.
Then:
\begin{itemize}
\item[(1)]
$w_\nu = 1$ for all $\nu \leq S$, or
\item[(2)]
$|x_{ik}|=|x_{il}|$ for all $i \in I$ and all $k, l \in \supp x_i$
\end{itemize}
\end{th}

Theorems~\ref{blocko} and \ref{blockl} show that in Orlicz and Lorentz
sequence spaces which are sufficiently different from $\ell_p$, a subspace
spanned by a block basis $\{f_j\}_j$ is $1$-complemented if and only if all
elements of the block basis $\{f_j\}_j$ have constant coefficients, i.e. if
they  satisfy Theorem~\ref{block}, about the most obvious form of
$1$-complemented subspaces in symmetric spaces.  We believe that in fact
Theorem~\ref{block} provides not only a sufficient condition but also a
necessary condition for a form of a $1$-complemented subspaces in sequence
spaces sufficiently different from $\ell_p$.

\subsection{Preservation of approximation properties by norm one projections}
\lb{ap}

In this final section we  mention briefly  problems analogous
to Problems~\ref{basis} and \ref{comp}, but concerned with structures which are
analogous to but weaker than unconditional bases or bases, namely finite
dimensional decompositions and approximation properties.

The approximation property appeared already in Banach's book \cite{Ban}
and it plays a fundamental role in the structure theory of Banach spaces,
see \cite{LT1} and the recent \cite{Cas} for the exposition of the
development and open problems in the theory. Here we just recall
the definitions essential for the statement of results below.

\begin{defn} 
\cite{CK90,GK97} (see also \cite[Sections~1.e and 1.g]{LT1} and \cite{Cas})
Let $X$ be a separable Banach space. We say that $X$ has the
{\it approximation property} (AP for short) if there is a net of finite
rank operators $T_\al$ so that $T_\al x \lra x$ for $x\in X$, uniformly
on compact sets. We say that $X$ has the
{\it bounded approximation property} (BAP for short) if this net can be
replaced by a sequence $T_n$, alternatively X has BAP 
if there is a sequence of finite
rank operators $T_n$ so that $T_n x \lra x$ for $x\in X$ and
$\sup_n\|T_n\|<\infty$.
A sequence $T_n$ with these properties is called an {\it approximating 
sequence}.
If $X$ has an approximating 
sequence $T_n$ with $\lim_{n\to\infty}\|T_n\|=1$ then we say that
$X$ has the
{\it metric approximation property} (MAP for short).
If $X$ has an approximating 
sequence $T_n$ with $\lim_{n\to\infty}\|I-2T_n\|=1$ then we say that
$X$ has the
{\it unconditional metric approximation property} (UMAP for short)
(here $I$ denotes the identity operator on $X$).
If $X$ has an approximating 
sequence $T_n$ with $T_mT_n=T_{\min(m,n)}$ then we say that
$X$ has a {\it finite dimensional decomposition} (FDD for short).
Let $1\le \la<\infty$. An FDD is called {\it $\la-$unconditional}
if $\|\sum_{n=1}^\infty \theta_n(T_{n+1}-T_n)\|\le \la$ for every
choice of signs $(\theta_n)_{n=1}^\infty$.
\end{defn} 

As mentioned above we refer the reader to \cite{Cas} for the analysis
of these interesting properties. Here we just want to note that, clearly,
spaces with bases have AP and FDD, spaces with monotone bases have MAP
and spaces with 1-unconditional bases have UMAP. It is a very nontrivial fact
that the reverse implications do not hold, see \cite{LT1,Cas}.

It is very natural to consider the following analogue of Problems~\ref{basis}
and \ref{comp}:

\begin{problem}
\lb{compap}
Let $X$ be a separable Banach space with one of the approximation
properties: AP, BAP, MAP, UMAP or FDD or a basis. Suppose that $Y$ is
complemented in $X$. Which of the approximation
properties must $Y$ satisfy?
\end{problem}

We note here that given any pair of Banach spaces $X, Y$ such that
$Y\subset X$ and $Y$ is complemented in $X$, it is possible to introduce
a new equivalent norm on $X$ so that $Y$ is 1-complemented in $X$ with the
new norm.
Thus Problem~\ref{compap} can be restated for 1-complemented subspaces $Y$ of
 $X$ without loosing the isomorphic nature of the problem.
We refer the reader to \cite{Cas} for an interesting account of what is 
known about Problem~\ref{compap}. Here we just quote one sample result
related to this problem:

\begin{th} \cite[Theorem 3.13]{Cas} (due to Pe\l czy\'nski \cite{P71}
and Johnson, Rosenthal, Zippin \cite{JRZ71})
A separable Banach space $Y$ has BAP \wtw $Y$ is isomorphic to a
complemented (equivalently, 1-complemented) subspace of a space with
a basis.
\end{th}

Isometric versions of Problem~\ref{compap} will necessarily deal with
isometric versions of approximation
properties i.e. MAP, UMAP, 1-unconditional FDD, monotone basis,
 1-unconditional  basis, similarly as Problems~\ref{monotone} and
\ref{1-unc}. 

Clearly 1-complemented  subspaces of  spaces with MAP have MAP.
Also Godefroy and Kalton \cite{GK97} observed that 
1-complemented  subspaces of  spaces with UMAP have UMAP.
Moreover they proved:

\begin{th} \cite[Corollary IV.4]{GK97}
Let $X$ be a separable Banach space. Then $X$ has UMAP 
\wtw for every $\e>0,$ $X$ is isometric to a 
1-complemented subspace of a space $V_\e$ with a $(1+\e)-$unconditional
FDD.
\end{th}

As far as we know, all other isometric versions of Problem~\ref{compap}
are open,
see also Section~\ref{ecp}.

We have for example:

\begin{problem}
 Does every 1-complemented subspace of a space 
 with a 1-unconditional
FDD have a 1-unconditional FDD (or any  unconditional FDD)? 
\end{problem} 

But note that Read \cite{Read} gave examples of spaces without FDD
which are complemented (equivalently, 1-complemented) in a space with an
FDD, or even in a space with a basis (cf. \cite{Cas}).


\end{document}